\documentclass[11pt]{article}
\usepackage{tikz}
\usepackage[all]{xy}
\usepackage{amsfonts}
\usepackage{amsfonts}
\usepackage{amsfonts}
\usepackage{mathrsfs}
\usepackage{bm}
\usepackage{cite}
\usepackage[colorlinks,linkcolor=blue,citecolor=blue]{hyperref}
\usepackage{amssymb,amsmath} % font used for R in Real numbers
\usepackage{makecell}
\usepackage{tikz}
 \allowdisplaybreaks
\usetikzlibrary{arrows,shapes,chains}
 \allowdisplaybreaks

\catcode`!=11
\let\!int\int \def\int{\displaystyle\!int}
\let\!lim\lim \def\lim{\displaystyle\!lim}
\let\!sum\sum \def\sum{\displaystyle\!sum}
\let\!sup\sup \def\sup{\displaystyle\!sup}
\let\!inf\inf \def\inf{\displaystyle\!inf}
\let\!cap\cap \def\cap{\displaystyle\!cap}
\let\!max\max \def\max{\displaystyle\!max}
\let\!min\min \def\min{\displaystyle\!min}
\let\!frac\frac \def\frac{\displaystyle\!frac}
\catcode`!=12

\let\oldsection\section
\renewcommand\section{\setcounter{equation}{0}\oldsection}

\allowdisplaybreaks
\def\pf{\it{Proof.}\rm\quad}

\def\N{\mathbb{N}}

\def\z{\zeta}

\def\L{{\rm Li}}

\newtheorem{defn}{Definition}[section]
\newtheorem{thm}{Theorem}[section]

\newtheorem{cor}[thm]{Corollary}
\newtheorem{con}[thm]{Conjecture}
\newtheorem{re}[thm]{Remark}

\DeclareMathOperator*{\cat}{\mathbf{Cat}}
\begin{document}
%%%%%%%%%%%%%%%%%%%% title %%%%%%%%%%%%%%%%%%%%%%%%%%%%%%%%%%%%%%%%%%%%%%%%
\title {\bf Explicit Relations between Multiple Zeta Values and Related Variants}
\author{
{Ce Xu \thanks{Email: cexu2020@ahnu.edu.cn}}\\[1mm]
 School of Mathematics and Statistics, Anhui Normal University,\\ \small Wuhu 241000, P.R. China\\
[5mm]
Dedicated to professor Masanobu Kaneko on the occasion of his 60th birthday}

\date{}
\maketitle \noindent{\bf Abstract} In this paper we present some new identities for multiple polylogarithms (abbr. MPLs) and multiple harmonic star sums (abbr. MHSSs) by using the methods
of iterated integral computations of logarithm functions. Then, by applying these formulas obtained, we establish some explicit relations between Kaneko-Yamamoto type multiple zeta values (abbr. K-Y MZVs), multiple zeta values (abbr. MZVs) and MPLs. Further, we find some explicit relations between MZVs and multiple zeta star values (abbr. MZSVs). Furthermore, we define an Ap\'{e}ry-type variant of MZSVs $\zeta^\star_B({\bf k})$ (called multiple zeta $B$-star values, abbr. MZBSVs) which involve MHSSs and central binomial coefficients, and establish some explicit connections among MZVs, alternating MZVs and MZBSVs by using the method of iterated integrals. Finally, some interesting consequences and illustrative examples are presented.
\\[2mm]


\section{Introduction and Notations}
We begin with some basic notation. A finite sequence ${\bf k} = (k_1,\ldots, k_r)$ of positive integers is called an \emph{index}. As usual, we
put
\[|{\bf k}|:=k_1+\cdots+k_r,\quad d({\bf k}):=r,\]
and call them the weight and the depth of ${\bf k}$, respectively. If $k_1>1$, ${\bf k}$ is called \emph{admissible}.

For an index ${\bf k}=(k_1,\ldots,k_r)$ and positive integer $n$, the \emph{multiple harmonic sums} (MHSs for short) and \emph{multiple harmonic star sums} (MHSSs for short) are defined by
\begin{align}
&{\zeta _n}( {\bf k})\equiv {\zeta _n}(k_1,k_2, \cdots ,k_r): = \sum\limits_{n \ge {n_1} > {n_2} >  \cdots  > {n_r} \ge 1} {\frac{1}{{n_1^{k_1}n_2^{k_2} \cdots n_r^{k_r}}}} ,\label{1.1}\\
&{\zeta_n ^ \star }({\bf k} )\equiv{\zeta_n ^ \star }(k_1,k_2,\ldots,k_r): = \sum\limits_{n \ge {n_1} \ge {n_2} \ge  \cdots  \ge {n_r} \ge 1} {\frac{1}{{n_1^{k_1}n_2^{k_2} \cdots n_r^{k_r}}}},\label{1.2}
\end{align}
when $n<k$, then ${\zeta_n}({\bf k}):=0$, and ${\zeta _n}(\emptyset )={\zeta^\star _n}(\emptyset ):=1$. When taking the limit $n\rightarrow \infty$ in (\ref{1.1}) and (\ref{1.2}), we get the so-called the \emph{multiple zeta values} (MZVs for short) and the \emph{multiple zeta star values} (MZSVs for short), respectively
\begin{align}
&{\zeta}( {\bf k}):=\lim_{n\rightarrow \infty}{\zeta _n}( {\bf k}),\label{1.3}\\
&{\zeta^\star}( {\bf k}):=\lim_{n\rightarrow \infty}{\zeta^\star_n}( {\bf k}),\label{1.4}
\end{align}
defined for an admissible index  ${\bf k}$ to ensure convergence of the series. The study of multiple zeta values began in the early 1990s with the works of Hoffman \cite{H1992} and Zagier \cite{DZ1994}. For an admissible index  ${\bf k}$, Hoffman \cite{H1992} called (\ref{1.3}) multiple harmonic series. Zagier \cite{DZ1994} called (\ref{1.3}) multiple zeta values since for $r=1$ they generalize the usual Riemann zeta values $\zeta(k)$. It has been attracted a lot of research on them in the last three decades (see, for example, the book of Zhao \cite{Z2016}).

For a non-empty index ${\bf k}=(k_1,\ldots,k_r)$, we write ${\bf k}^\star$
for the formal sum of $2^{r-1}$ indices of the form
$(k_1\bigcirc \cdots \bigcirc k_r)$, where each $\bigcirc$ is
replaced by `\,,\,' or `+'.
We also put $\varnothing^\star=\varnothing$. Then, we have
$\zeta_n^\star({\bf k})=\zeta_n({\bf k^\star)}$ for ${\bf k}\in \N^r$.

Recently, Kaneko and Yamamoto \cite{KY2018} introduced and studied a new kind of multiple zeta values
\begin{align}\label{1.5}
\zeta({\bf k}\circledast{\bf l}^\star)&=\sum\limits_{0<m_r<\cdots<m_1=n_1\geq \cdots \geq n_s>0} \frac{1}{m_1^{k_1}\cdots m_r^{k_r}n_1^{l_1}\cdots n_s^{l_s}}\nonumber\\
&=\sum\limits_{n=1}^\infty \frac{\zeta_{n-1}(k_2,\ldots,k_r)\zeta^\star_n(l_2,\ldots,l_s)}{n^{k_1+l_1}},
\end{align}
where ${\bf k}=(k_1,\ldots,k_r)$ and ${\bf l}=(l_1,\ldots,l_s)$ are any arrays of positive integers. We call them \emph{K-Y multiple zeta values} (K-Y MZVs for short). Note that we used the opposite convention of the original definition of K-Y MZVs. It is clear that the left-hand side is a certain integral which can be written as a $\mathbb{Z}$-linear combination of MZVs.
The same holds if we use the series expression on the right. Kaneko and Yamamoto \cite{KY2018} presented a
new ``integral=series" type identity of multiple zeta values, and conjectured that this identity is
enough to describe all linear relations of multiple zeta values over  $\mathbb{Q}$.
Here ``${\bf k}\circledast{\bf l}$" defined by the `\emph{circle harmonic shuffle product}' of indexes ${\bf k}=(k_1,\ldots,k_r)$ and ${\bf l}=(l_1,\ldots,l_s)$:
\begin{align*}
&{\bf k}\circledast{\bf l}:=
\bigl(k_1+l_1,(k_2,\ldots,k_{r})*(l_2,\ldots,l_{s})\bigr),
\end{align*}
and ``${\bf k}\ast{\bf l}$" defined by the `\emph{harmonic shuffle product}' of indexes ${\bf k}=(k_1,\ldots,k_r)$ and ${\bf l}=(l_1,\ldots,l_s)$:
\begin{align*}
&\varnothing *\bf k=\bf k*\varnothing=\bf k,\\
&{\bf k*\bf l}
=\bigl(k_1,(k_2,\ldots,k_{r})*(l_1,\ldots,l_s)\bigr)\\
&\quad\quad\quad+\bigl(l_1,(k_1,\ldots,k_r)*(l_2,\ldots,l_{s})\bigr)\\
&\quad\quad\quad+\bigl(k_1+l_1,(k_2,\ldots,k_{r})*(l_2,\ldots,l_{s})\bigr),
\end{align*}
where $\emptyset$ denotes the unique index of depth $0$.

For convenience, by ${\left\{ {{s_1}, \ldots ,{s_j}} \right\}_d}$ we denote the sequence of depth $dj$ with $d$ repetitions of ${\left\{ {{s_1}, \ldots ,{s_j}} \right\}}$. For example,
\begin{align*}
{\left\{s_1,s_2,s_3 \right\}_4}=\left\{s_1,s_2,s_3,s_1,s_2,s_3,s_1,s_2,s_3,s_1,s_2,s_3\right\}.
\end{align*}

Clearly, setting ${\bf k}=(1)$ or ${\bf l}=(1)$ yield
\begin{align*}
&\zeta({\bf k}\circledast(1)^\star)=\zeta(k_1+1,k_2,\ldots,k_r),\\
&\zeta((1)\circledast{\bf l}^\star)=\zeta^\star(l_1+1,l_2,\ldots,l_s).
\end{align*}
Specially, the \emph{Arakawa-Kaneko zeta value}
\begin{align*}
\xi(p;{\bf k})=\sum\limits_{n=1}^\infty \frac{\zeta_{n-1}(k_2,\ldots,k_r)\zeta^\star_n(\{1\}_{p-1})}{n^{k_1+1}}=\zeta({\bf k}\circledast \{\underbrace{1,\ldots,1}_{p}\}^\star)
\end{align*}
is also a special case of K-Y MZV (see \cite{Ku2010}), where $p,k_1,\ldots,k_r\in \N$. Here the \emph{Arakawa-Kaneko function} is defined, for $\Re(s)>0$ and positive integers $k_1,k_2,...,k_r\ (r\in\N)$, by (\cite{AM1999})
\begin{align}\label{1.6}
\xi(s;k_1,k_2\ldots,k_r):=\frac{1}{\Gamma(s)} \int\limits_{0}^\infty \frac{t^{s-1}}{e^t-1}{\rm Li}_{k_1,k_2,\ldots,k_r}(1-e^{-t})dt,
\end{align}
where the \emph{multiple polylogarithm} ${\rm Li}_{k_1,k_2,\ldots,k_r}(x)$ (MPL for short) is defined by ($(k_1,x)\neq (1,1)$)
\begin{align}\label{1.7}
{\mathrm{Li}}_{k_1,k_2,\ldots,k_r}\left(x\right): = \sum\limits_{n_1>n_2>\cdots>n_r\geq 1} {\frac{x^{n_1}}{{n_1^{k_1}n_2^{k_2} \cdots n_r^{k_r}}}}  \quad x\in [-1,1].
\end{align}
Some related results for Arakawa-Kaneko functions and related functions may be seen in the works of \cite{BH2011,ChenKW19,CC2010,CC2015,I2016,KTA2018,Ku2010,Yo2014,Yo2015} and references therein.

The primary goals of this paper are to study the explicit relations of MZVs and their related variants, such as K-Y MZVs $\zeta({\bf k}\circledast{\bf l}^\star)$. Then using these explicit relations, we establish some explicit formulas of MZVs and MZSVs.

The remainder of this paper is organized as follows.

In Section \ref{sec2}, we establish several iterated integral formulas of multiple polylogarithms and define a parametric variation of multiple harmonic star sums.

In Section \ref{sec3}, we will establish general relations between K-Y MZVs and MZVs by using the method of iterated integration. In particular, we show that
\begin{align*}
\zeta((\{1\}_k)\circledast(1,\{\{1\}_{m-1},2\}_{p-1},\{1\}_m)^\star)\in\mathbb{Q}[\text{Riemann Zeta Values}]
\end{align*}
and
\[\zeta((\{1\}_m,2,\{1\}_b)\circledast (1,\{\{1\}_{m-1},2\}_{p-1},\{1\}_m)^\star)\in\mathbb{Q}[\text{Riemann Zeta Values}].\]
Moreover,  some interesting consequences and illustrative examples are considered.

In Section \ref{sec4}, for an index ${\bf k}$ we define the Ap\'{e}ry-type variant of MZSV
\begin{align}
\zeta^\star_B({\bf k})\equiv\zeta^\star_B(k_1,k_2,\ldots,k_r):=\sum\limits_{n=1}^\infty \frac{\zeta_n^\star(k_2,\ldots,k_r)}{n^{k_1}4^n}\binom{2n}{n},
\end{align}
we call them \emph{multiple zeta $B$-star values} (MZBSVs for short). Using a similar method to Section \ref{sec3}, we will prove several explicit relations between $\zeta^\star_B({\bf k})$ and alternating MZ(S)Vs. For non-zero integers $k_j\ (j=1,2,\ldots,r)$, the \emph{alternating multiple zeta values} (AMZVs for short)  and \emph{alternating multiple zeta star values} (AMZSVs for short) are defined by
\begin{align}
&\zeta(k_1, \ldots, k_r):=\sum\limits_{n_1>\cdots>n_r\geq 1}\prod\limits_{j=1}^r n_j^{-|k_j|}{\rm sgn}(k_j)^{n_j},\\
&\zeta^\star(k_1, \ldots, k_r):=\sum\limits_{n_1\geq\cdots\geq n_r\geq 1}\prod\limits_{j=1}^r n_j^{-|k_j|}{\rm sgn}(k_j)^{n_j},
\end{align}
where for convergence $|k_1|+\cdots+|k_j|> j$ for $j= 1, 2, \ldots, r$ (we allow $k_1=-1$), and
\[{\rm sgn}(k_j):=\begin{cases}
   1  & \text{\;if\;} k_j>0,  \\
   -1, & \text{\;if\;} k_j<0.
\end{cases}\]
When ${\rm sgn}(s_j)=-1$,  by placing a bar over the
corresponding integer exponent $|k_j|$. For example, we write
\begin{align*}
{\zeta}( {\bar 2,3,\bar 1,4} )={\zeta}( {-2,3,- 1,4})\quad{\rm and}\quad \z^\star( {\bar 3,2,\bar 4,1})=\z^\star( {-3,2,-4,1}).
\end{align*}

\section{Iterated integrals of Multiple Polylogarithms and Parametric Multiple Harmonic Star Sums}\label{sec2}

In this section, we establish several iterated integral identities of multiple polylogarithms and parametric multiple harmonic star sums.
We begin with some basic notations.
For any index ${\bf m}:=(m_1,m_2,\ldots,m_p)$, we define \begin{align*}
&{\bf m}_j:=(m_1,m_2+1,\ldots,m_j+1),\ {\bf m}_1:=(m_1),\\
&\overrightarrow{{\bf m}}_{i,j}:=(m_i,m_{i+1},\ldots,m_j),\\
&\overleftarrow{{\bf m}}_{i,j}:=(m_j,m_{j-1},\ldots,m_i),
\end{align*}
and
\begin{align*}
|{\bf m}|_j:=m_1+m_2+\cdots+m_j,
\end{align*}
where $1\leq i\leq j\leq p$. If $i>j$, then we let $\overleftarrow{{\bf m}}_{i,j}:=\emptyset$. In definition of ${\bf m}_j$, we allow $m_i=0\ (i=2,3,\ldots,j)$.

Moreover, for two index ${\bf m}:=(m_1,m_2,\ldots,m_p)$ and ${\bf n}:=(n_1,n_2,\ldots,n_p)$, we let
\begin{align*}
&\overrightarrow{({\bf m+n})}_{i,j}:=(m_i+n_i,\ldots,m_{j-1}+n_{j-1},m_j+n_j),\\
&\overleftarrow{({\bf m+n})}_{i,j}:=(m_j+n_j,m_{j-1}+n_{j-1},\ldots,m_i+n_i).
\end{align*}
In particular, if ${\bf n}=(\{1\}_p)$, then
\begin{align*}
&\overrightarrow{({\bf m+1})}_{i,j}:=(m_i+1,\ldots,m_{j-1}+1,m_j+1),\\
&\overleftarrow{({\bf m+1})}_{i,j}:=(m_j+1,m_{j-1}+1,\ldots,m_i+1).
\end{align*}

\subsection{Identities for Multiple Polylogarithms}

According to the definition of multiple polylogarithms, we have
\begin{align}\label{2.1}
\frac{d}{dx}{\mathrm{Li}}_{{{k_1},{k_2}, \cdots ,{k_r}}}(x)= \left\{ {\begin{array}{*{20}{c}} \frac{1}{x} {\mathrm{Li}}_{{{k_1-1},{k_2}, \cdots ,{k_r}}}(x)
   {,\ \ k_1>1,}  \\
   {\frac{1}{1-x}{\mathrm{Li}}_{{{k_2}, \cdots ,{k_r}}}(x),\;\;\;k_1 = 1.}  \\
\end{array} } \right.
\end{align}
Hence, applying (\ref{2.1}) we obtain the following iterated integral expression
\begin{align}\label{2.2}
{\mathrm{Li}}_{{{k_1},{k_2}, \cdots ,{k_r}}}\left( x \right)=\int\limits_{0}^x \underbrace{\frac{dt}{t}\cdots\frac{dt}{t}}_{k_1-1}\frac{dt}{1-t}\underbrace{\frac{dt}{t}\cdots\frac{dt}{t}}_{k_2-1}\frac{dt}{1-t}\cdots
\underbrace{\frac{dt}{t}\cdots\frac{dt}{t}}_{k_r-1}\frac{dt}{1-t},
\end{align}
where $0<x<1$ and
$$\int\limits_{0}^x f_1(t)dtf_2(t)dt\cdots f_r(t)dt:=\int\limits_{0<t_r<\cdots<t_1<x}f_1(t_1)f_2(t_2)\cdots f_r(t_r)dt_1dt_2\cdots dt_r.$$
Note the fact that
\begin{align*}
\int\limits_{x_1<t_r<\cdots<t_1<x_2} \frac{dt_1\cdots dt_r}{t_1\cdots t_r} =\frac{1}{r!}\log^r\left(\frac{x_2}{x_1}\right).
\end{align*}
Thus, (\ref{2.2}) can be rewritten in the form
\begin{align}\label{2.3}
{\mathrm{Li}}_{{{k_1},{k_2}, \cdots ,{k_r}}}\left( x \right)=\prod\limits_{j=1}^r\frac{1}{\Gamma(k_j)}\int\limits_{0<t_r<\cdots<t_1<t_0=x} \prod\limits_{j=1}^r\left\{\frac{\log^{k_j-1}\left(\frac{t_{j-1}}{t_j}\right)}{1-t_j}dt_j\right\}.
\end{align}
For convenience, we use the following notations
\begin{align*}
\Omega_{p}^{k}(x):=\frac{\log^{k-1}\left(\frac{x}{1-t_p}\right)}{t_p}\quad {\rm and}\quad \Omega_{p,p-1}^{k}:=\frac{\log^{k-1}\left(\frac{1-t_p}{1-t_{p-1}}\right)}{t_{p-1}},
\end{align*}
where $p$ and $k$ are positive integers, if $k=1$ we let $\Omega_{p}(x)\equiv\Omega^{1}_{p}(x)$ and $\Omega_{p,p-1}\equiv\Omega_{p,p-1}^{1}$. Let
\begin{align*}
E_{r}(x):=\{(t_1,\ldots,t_r)\mid x<t_r<\cdots<t_1<1\},\quad r\in\mathbb{N}.
\end{align*}
Therefore, in (\ref{2.3}), changing variable  $t_i\mapsto 1-t_{r+1-i}\ (i=1,2,\ldots,r)$ gives
\begin{align}\label{2.4}
{\mathrm{Li}}_{{{k_1},{k_2}, \cdots ,{k_r}}}\left( x \right)=\frac{1}{\prod\limits_{j=1}^r \Gamma(k_j) } \int\nolimits_{E_r(1-x)} \Omega_{r}^{k_1}(x)\Omega_{r,r-1}^{k_2}\cdots\Omega_{2,1}^{k_r}dt_1\cdots dt_r.
\end{align}
\begin{thm} For $k,r\in \N$, we have
\begin{align}\label{2.5}
{\rm Li}_{k,\{1\}_{r-1}}(x)=&\sum\limits_{j=1}^{k-1} \frac{\log^{k-1-j}(x)}{(k-1-j)!}\left\{\zeta(r+1,\{1\}_{j-1})-\sum\limits_{i=0}^{r-1}\frac{(-1)^i}{i!}\log^i(1-x){\rm Li}_{r+1-i,\{1\}_{j-1}}(1-x) \right\}\nonumber\\
&+(-1)^r \frac{\log^{k-1}(x)\log^r(1-x)}{(k-1)!r!}.
\end{align}
\end{thm}
\pf Setting $k_1=k,k_2=\cdots=k_r=1$ in (\ref{2.4}) yields
\begin{align}\label{2.6}
{\rm Li}_{k,\{1\}_{r-1}}(x)&=\frac{1}{(k-1)! } \int\nolimits_{D_{x,{\bf t}}} \Omega_{r}^{k}(x)\Omega_{r,r-1}^{1}\cdots\Omega_{2,1}^{1}dt_1\cdots dt_r\nonumber\\
&=\frac{1}{(k-1)!}\sum\limits_{j=0}^{k-1} \binom{k-1}{j}\log^{k-1-j}(x) \int\nolimits_{D_{x,{\bf t}}} \Omega_{r}^{j+1}(1)\Omega_{r,r-1}^{1}\cdots\Omega_{2,1}^{1}dt_1\cdots dt_r.
\end{align}
By direct calculations, we can find that
\begin{align}\label{2.7}
{\rm Li}_{\{1\}_r}(x)=\int\nolimits_{D_{x,{\bf t}}} \Omega_{r}^{1}(1)\Omega_{r,r-1}^{1}\cdots\Omega_{2,1}^{1}dt_1\cdots dt_r=\frac{(-1)^r}{r!}\log^r(1-x)
\end{align}
and $(j\geq 1)$
\begin{align}\label{2.8}
&\int\nolimits_{D_{x,{\bf t}}} \Omega_{r}^{j+1}(1)\Omega_{r,r-1}^{1}\cdots\Omega_{2,1}^{1}dt_1\cdots dt_r\nonumber\\=&j!\zeta(r+1,\{1\}_{j-1})-j!\sum\limits_{i=0}^{r-1}\frac{(-1)^i}{i!}\log^i(1-x){\rm Li}_{r+1-i,\{1\}_{j-1}}(1-x).
\end{align}
Thus, substituting (\ref{2.7}) and (\ref{2.8}) into (\ref{2.6}), we immediately obtain the formula (\ref{2.5}).\hfill$\square$

\begin{re} This theorem generalizes \cite[Theorem 8]{AM1999}, where the corresponding formula for ${\rm Li}_{k,\{1\}_{r-1}}(x)$ is
\begin{align}\label{2.9}
{\rm Li}_{k,\{1\}_{r-1}}(x)=&(-1)^{k-1} \sum\limits_{\scriptstyle j_1+\cdots+j_k=r+k, \hfill \atop
  \scriptstyle\ \  j_1,\ldots,j_k\geq 1 \hfill} {\rm Li}_{\{1\}_{j_k-1}}(x){\rm Li}_{j_1,\ldots,j_{k-1}}(1-x)\nonumber\\
  &+\sum\limits_{j=0}^{k-2} (-1)^j \zeta(k-j,\{1\}_{r-1}){\rm Li}_{\{1\}_j}(1-x).
\end{align}
Comparing (\ref{2.5}) with (\ref{2.9}), and noting that ${\rm Li}_{\{1\}_j}(x)=\frac{(-1)^j}{j!}\log^j(1-x)$, we obtain
\begin{align}\label{2.10}
(-1)^k \sum\limits_{\scriptstyle j_1+\cdots+j_{k-1}=r+k, \hfill \atop
  \scriptstyle\ \  j_1,\ldots,j_{k-1}\geq 1 \hfill} {\rm Li}_{j_1,\ldots,j_{k-1}}(1-x)=\sum\limits_{j=1}^{k-1} \frac{\log^{k-1-j}(x){\rm Li}_{r+2,\{1\}_{j-1}}(1-x)}{(k-1-j)!},
\end{align}
where $r\in \N_0:=\N\cup \{0\}$ and $k\in\N$.
\end{re}

Recently, Kaneko and Tsumura \cite{KT2018} gave the following more general result.
\begin{thm}\label{thm2.3}(Kaneko-Tsumura \cite{KT2018}) Let ${\bf k}$ be any index. Then we have
\begin{align}\label{2.11}
{\rm Li}_{\bf k}(1-x)=\sum\limits_{{\bf k},j\geq 0} c_{\bf k}({\bf k'};j){\rm Li}_{\{1\}_j}(1-x){\rm Li}_{\bf k'}(x),
\end{align}
where the sums on the right runs over indices ${\bf k'}$ and integers $j\geq 0$ that satisfy $|{\bf k'}|+j\leq |{\bf k}|$, and $c_{\bf k}({\bf k'};j)$ is a $\mathbb{Q}$-linear combination of multiple zeta values of weight $|{\bf k}|- |{\bf k'}|-j$. We understand ${\rm Li}_{\emptyset}(x)=1$ and $|\emptyset|=0$ for the empty index $\emptyset$, and the constant $1$ is regarded as a multiple zeta value of weight $0$.
\end{thm}

In particular, the author and Pallewattaa \cite[Thm. 2.9]{PX2019} gave the following an explicit formula
\begin{align}\label{eb1}
{\rm Li}_{\{1\}_a,2,\{1\}_b}(1-x)&=\sum_{j=0}^a (-1)^j \binom{j+b+1}{j}\z(j+b+2){\rm Li}_{\{1\}_{a-j}}(1-x)\nonumber\\
&\quad-(-1)^a \sum_{j=0}^{b+1} \binom{j+a}{j} {\rm Li}_{\{1\}_{b+1-j}}(1-x){\rm Li}_{a+1+j}(x),
\end{align}
where $a,b\in\N_0$.

Next, we use the identities (\ref{2.5}), (\ref{2.11}) and (\ref{eb1}) to establish some iterated integral formulas involving multiple polylogarithms.
\begin{thm}\label{thm2.4} For ${\bf m}=(m_1,\ldots,m_p)\in \N_0^p$ and ${\bf k}=(k_1,k_2,\ldots,k_r)\in \N^r$,
\begin{align}\label{2.12}
&\int\nolimits_{E_p(x)} \Omega_p^{m_p+1}(1-x)\Omega_{p,p-1}^{m_{p-1}+1}\cdots \Omega_{2,1}^{m_1+1}{\rm Li}_{\bf k}(t_1)  dt_1\cdots dt_p\nonumber\\
&=m_1!\cdots m_p!\sum\limits_{{\bf k'},j\geq 0}c_{\bf k}({\bf k'};j)\sum\limits_{j_0+\cdots+j_p=j,\atop j_0,\ldots,j_p \geq 0} \frac{(-1)^{j_0}}{j_0!}\left\{\prod\limits_{l=1}^p \binom{m_l+j_l}{j_l}\right\}\nonumber\\&\quad\quad\quad\quad\quad\quad\quad\quad\times\log^{j_0}(1-x){\rm Li}_{\overleftarrow{({\bf m+j+1})}_{1,p},{\bf k'}}(1-x),
\end{align}
Here, the the sums is over indices ${\bf k'}$ and integers $j\geq 0$ that satisfy $|{\bf k'}|+j\leq |{\bf k}|$, and $c_{\bf k}({\bf k'};j)$ is a $\mathbb{Q}$-linear combination of multiple zeta values of weight $|{\bf k}|- |{\bf k'}|-j$.
\end{thm}
\pf Applying the change of variables $t_i \rightarrow 1-(1-x)t_1\cdots t_{p+1-i}$, and using (\ref{2.11}) then the integral on the left hand sides of (\ref{2.12}) can be rewritten as
\begin{align}\label{2.13}
&\int\nolimits_{E_p(x)} \Omega_p^{m_p+1}(1-x)\Omega_{p,p-1}^{m_{p-1}+1}\cdots \Omega_{2,1}^{m_1+1}{\rm Li}_{\bf k}(t_1)  dt_1\cdots dt_p\nonumber\\
&=(-1)^{m_1+\cdots+m_p} \underbrace{\int\limits_{0}^1\cdots\int\limits_{0}^1}_{p} \left\{\prod\limits_{j=1}^p\frac{\log^{m_{p+1-j}}(t_j)}{\left(((1-x)t_1\cdots t_{j-1})^{-1}-t_j\right)}\right\}{\rm Li}_{\bf k}(1-(1-x)t_1\cdots t_p) dt_1\cdots dt_p\nonumber\\
&=(-1)^{m_1+\cdots+m_p} \sum\limits_{{\bf k},j\geq 0} c_{\bf k}({\bf k'};j)\frac{(-1)^j}{j!}\underbrace{\int\limits_{0}^1\cdots\int\limits_{0}^1}_{p} \left\{\prod\limits_{i=1}^p\frac{\log^{m_{p+1-i}}(t_i)}{\left(((1-x)t_1\cdots t_{i-1})^{-1}-t_i\right)}\right\}\nonumber\\&\quad\quad\quad\quad\quad\quad\quad\quad\quad\quad\quad\quad\quad\quad\quad\quad\times{\log^j((1-x)t_1\cdots t_p)}{\rm Li}_{\bf k'}((1-x)t_1\cdots t_p) dt_1\cdots dt_p\nonumber\\
&=(-1)^{m_1+\cdots+m_p} \sum\limits_{{\bf k},j\geq 0} c_{\bf k}({\bf k'};j)(-1)^j\sum\limits_{j_0+\cdots+j_p=j,\atop j_0,\ldots,j_p \geq 0}\frac{\log^{j_0}(1-x)}{j_0!\cdots j_p!} \nonumber\\ &\quad\quad\quad\quad\times \underbrace{\int\limits_{0}^1\cdots\int\limits_{0}^1}_{p} \left\{\prod\limits_{i=1}^p\frac{\log^{m_{p+1-i}+j_i}(t_i)}{\left(((1-x)t_1\cdots t_{i-1})^{-1}-t_i\right)}\right\}{\rm Li}_{\bf k'}((1-x)t_1\cdots t_p) dt_1\cdots dt_p.
\end{align}
From (\ref{2.1}), we find that
\begin{align*}
\frac{(-1)^{k_1-1}}{\Gamma(k_1)}\int\limits_{0}^1x \frac{\log^{k_1-1}(t){\rm Li}_{k_2,\ldots,k_r}(xt)}{1-xt}dt={\mathrm{Li}}_{{{k_1},{k_2}, \cdots ,{k_r}}}\left( x \right).
\end{align*}
Further, we have
\begin{align}\label{2.14}
{\rm Li}_{\bf k,m}(x)=\prod\limits_{j=1}^r\frac{(-1)^{k_j-1}}{\Gamma(k_j)}\underbrace{\int\limits_{0}^1\cdots\int\limits_{0}^1}_{r} \left\{\prod\limits_{j=1}^r\frac{\log^{k_j-1}(t_j)}{\left(x\prod\limits_{l=1}^{j-1}t_l\right)^{-1}-t_j}\right\}{\rm Li}_{\bf m}\left(x\prod\limits_{j=1}^r t_j\right)\prod\limits_{j=1}^r dt_j,
\end{align}
where $t_0:=1$. Hence, combining (\ref{2.13}) with (\ref{2.14}), we obtain the desired formula (\ref{2.12}).\hfill$\square$

Letting ${\bf k}=(k,\{1\}_{r-1})$ in (\ref{2.12}) and using (\ref{2.5}), we can get the following corollary.
\begin{cor} For ${\bf m}=(m_1,\ldots,m_p)\in \N_0^p$ and $k,r\in \N$,
\begin{align}\label{2.15}
&\int\nolimits_{E_p(x)} \Omega_p^{m_p+1}(1-x)\Omega_{p,p-1}^{m_{p-1}+1}\cdots \Omega_{2,1}^{m_1+1}{\rm Li}_{k,\{1\}_{r-1}}(t_1)  dt_1\cdots dt_p\nonumber\\
&=m_1!\cdots m_p!(-1)^{k-1}\sum\limits_{j_1+\cdots+j_k=r,\atop j_1,\ldots,j_k\geq 0}\sum\limits_{i_0+\cdots+i_p=j_k,\atop i_0,\ldots,i_p \geq 0} \frac{(-1)^{i_0}}{i_0!}\left\{\prod\limits_{l=1}^p \binom{m_l+i_l}{i_l}\right\}\log^{i_0}(1-x)\nonumber\\&\quad\quad\quad\quad\quad\quad\quad\quad\quad\quad\quad\quad\times{\rm Li}_{\overleftarrow{({\bf m+i+1})}_{1,p},\overrightarrow{({\bf j+1})}_{1,k-1}}(1-x)\nonumber\\&\quad+m_1!\cdots m_p!\sum\limits_{j=0}^{k-2} (-1)^j \zeta(k-j,\{1\}_{r-1}){\rm Li}_{\overleftarrow{({\bf m+1})}_{1,p},\{1\}_j}(1-x),
\end{align}
where ${\bf i}:=(i_1,\ldots,i_p)$ and ${\bf j}:=(j_1,\ldots,j_{k-1})$.
\end{cor}

Setting ${\bf k}=(\{1\}_a,2,\{1\}_b)$ in (\ref{2.12}) and applying (\ref{eb1}) yields the following corollary.
\begin{cor} For ${\bf m}=(m_1,\ldots,m_p)\in \N_0^p$ and $a,b\in \N_0$,
\begin{align}\label{eb2}
&\int\nolimits_{E_p(x)} \Omega_p^{m_p+1}(1-x)\Omega_{p,p-1}^{m_{p-1}+1}\cdots \Omega_{2,1}^{m_1+1}{\rm Li}_{\{1\}_a,2,\{1\}_b}(t_1)  dt_1\cdots dt_p\nonumber\\
&=m_1!\cdots m_p! \sum_{j=0}^a (-1)^j \binom{j+b+1}{j}\z(j+b+2)\sum_{i_0+\cdots+i_p=a-j,\atop i_0,\ldots,i_p\geq 0} \frac{(-1)^{i_0}}{i_0!}\left\{\prod\limits_{l=1}^p \binom{m_l+i_l}{i_l}\right\}\nonumber\\
&\quad\quad\quad\quad\quad\quad\quad\quad\quad\quad\quad\quad\quad\quad\quad\quad\quad\quad\times\log^{i_0}(1-x){\rm Li}_{\overleftarrow{({\bf m+i+1})}_{1,p}}(1-x)\nonumber\\
&\quad-m_1!\cdots m_p!(-1)^a \sum_{j=0}^{b+1} \binom{j+a}{j}\sum_{i_0+\cdots+i_p=b+1-j,\atop i_0,\ldots,i_p\geq 0} \frac{(-1)^{i_0}}{i_0!}\left\{\prod\limits_{l=1}^p \binom{m_l+i_l}{i_l}\right\}\nonumber\\
&\quad\quad\quad\quad\quad\quad\quad\quad\quad\quad\quad\quad\quad\quad\quad\quad\times\log^{i_0}(1-x){\rm Li}_{\overleftarrow{({\bf m+i+1})}_{1,p},a+1+j}(1-x),
\end{align}
where ${\bf i}:=(i_1,\ldots,i_p)$.
\end{cor}

\subsection{Identities for Parametric Multiple Harmonic Star Sums}
For positive integers $m_1,\ldots,m_p$ and real $x\in [-1,1]$, we define a \emph{parametric multiple harmonic star sum} (PMHSS for short) $\zeta^\star_n(m_1,\cdots,m_{p-1},m_p;x)$ by
\begin{align*}
\z_n^\star(m_1,\cdots,m_{p-1},m_p;x):=\sum\limits_{n\geq n_1\geq \cdots \geq n_p\geq 1}\frac{x^{n_p}}{n^{m_1}_1\cdots n^{m_{p-1}}_{p-1}n^{m_p}_p},
\end{align*}
where $\zeta^\star_n(\emptyset;x):=x^n$.

It is clear that according to definition, by a direct calculation, we have the relations
\begin{align}\label{2.16}
\frac{d}{dx}\z_n^\star(m_1,\cdots,m_{p-1},m_p;x) = \left\{{\begin{array}{*{20}{c}}\ \ \ \quad\quad   \frac{1}{x} \z_n^\star(m_1,\cdots,m_{p-1},m_p;x)
   {,\ \ \ \quad\quad\quad m_p>1,}  \\
   {\frac{\z_n^\star(m_1,\cdots,m_{p-1})-\z_n^\star(m_1,\cdots,m_{p-1};x)}{1-x},\;\;\;m_p = 1.}  \\
\end{array} } \right.
\end{align}

\begin{defn}
For an index ${\bf m}=(m_1,\ldots,m_p)$, its Hoffman dual is
the index ${\bf m}^\vee=(m'_1,\ldots,m'_{p'})$ determined by
$|{\bf m}|:=m_1+\cdots+m_p=m'_1+\cdots+m'_{p'}$ and
\begin{equation*}
\{1,2,\ldots,|{\bf m}|-1\}
=\{m_1,m_1+m_2,\ldots,m_1+\cdots+m_{p-1}\}
\amalg\{m'_1,m'_1+m'_2,\ldots,m'_1+\cdots+m'_{p'-1}\}.
\end{equation*}
\end{defn}
For example, we have
\begin{align*}
({1,1,2,1})^\vee=(3,2)\quad\text{and}\quad ({1,2,1,1})^\vee=(2,3).
\end{align*}

\begin{thm} For positive integer $p$ and real $x\in (0,1]$,
\begin{align}\label{2.17}
\z_n^\star({\bf m}_p^v;x)=&\frac{(-1)^p n}{m_1!\cdots m_p!} \int\nolimits_{E_p(x)} \Omega_p^{m_p+1}(1-x)\Omega_{p,p-1}^{m_{p-1}+1}\cdots \Omega_{2,1}^{m_1+1}t_1^n dt_1\cdots dt_p\nonumber\\
&+ \sum\limits_{j=1}^{p} (-1)^{p-j}\z_n^\star({\bf m}_j^v){\rm Li}_{\overleftarrow{({\bf m+1})}_{j+1,p}}(1-x),
\end{align}
where ${\bf m}_p=(m_1,m_2+1,\ldots,m_p+1)$ and $m_1\geq 1,\ m_j\geq 0\ (j=2,3,\ldots,p)$.
\end{thm}
\pf The proof is by induction on $p$. Note the fact that the case $p=1$ is well known (see \cite{Xu2017}). For $p=2$, if $m_2=0$, by (\ref{2.16}), we have
\begin{align*}
\frac{d}{dx} \z_n^\star((m_1,1)^v;x)=\frac{\z_n^\star(m_1^v;x)}{x}
\end{align*}
and
\begin{align*}
&\frac{d}{dx} \left\{\frac{n}{m_1!} \int\nolimits_{E_2(x)} \Omega_2(1-x)\Omega_{2,1}^{m_1+1}t_1^n dt_1dt_2-\z(m_1^v)\L_1(1-x)+\z_n^\star((m_1,1)^v) \right\}\\
&=\frac{\z_n^\star(m_1^v;x)}{x}.
\end{align*}
Hence, we know that
\begin{align*}
\z_n^\star((m_1,1)^v;x)=\frac{n}{m_1!} \int\nolimits_{E_2(x)} \Omega_2(1-x)\Omega_{2,1}^{m_1+1}t_1^n dt_1dt_2-\z_n^\star(m_1^v)\L_1(1-x)+\z_n^\star((m_1,1)^v) +C_1,
\end{align*}
where $C_1$ is a constant. If letting $x\rightarrow 1$, then we get $C_1=0$.
Thus, the formula (\ref{2.17}) holds for $p=2$ and $m_2=0$. Similarly, if $m_2\geq 1$, we assume (\ref{2.17}) holds for $(m_1,m_2)$, then from (\ref{2.16}) we deduce
\begin{align*}
\frac{d}{dx} \z_n^\star((m_1,m_2+1)^v;x)=\frac{\z_n^\star((m_1,m_2)^v)-\z_n^\star((m_1,m_2)^v;x)}{1-x}
\end{align*}
and
\begin{align*}
&\frac{d}{dx} \left\{\frac{n}{m_1!m_2!} \int\nolimits_{E_2(x)} \Omega_2^{m_2+1}(1-x)\Omega_{2,1}^{m_1+1}t_1^n dt_1dt_2-\z_n^\star(m_1^v)\L_{m_2+1}(1-x)+\z_n^\star((m_1,m_2+1)^v) \right\}\\
&=\frac{\z_n^\star((m_1,m_2)^v)-\z_n^\star((m_1,m_2)^v;x)}{1-x}.
\end{align*}
Therefore, we obtain
\begin{align*}
\z_n^\star((m_1,m_2+1)^v;x)&=\frac{n}{m_1!m_2!} \int\nolimits_{E_2(x)} \Omega_2^{m_2+1}(1-x)\Omega_{2,1}^{m_1+1}t_1^n dt_1dt_2\\&\quad-\z_n^\star(m_1^v)\L_{m_2+1}(1-x)+\z_n^\star((m_1,m_2+1)^v)+C_2,
\end{align*}
where $C_2$ is a constant. Further, setting $x\rightarrow 1$ yields $C_2=0$. So, for $p=2$, the formula (\ref{2.17}) holds. Now, we assume the formula (\ref{2.17}) holds for index $(m_1,m_2+1,\ldots,m_{p-1}+1)^v$. If $m_p=0$, by direct calculations, we find that
\begin{align*}
&\frac{d}{dx} \z_n^\star(({\bf m}_{p-1},1)^v;x)=\frac{\z_n^\star((m_1,m_2+1,\ldots,m_{p-1}+1)^v;x)}{x}\\
&=\frac{d}{dx}\left\{\begin{array}{l} \frac{(-1)^p n}{m_1!\cdots m_{p-1}!} \int\nolimits_{E_p(x)} \Omega_p(1-x)\Omega_{p,p-1}^{m_{p-1}+1}\cdots \Omega_{2,1}^{m_1+1}t_1^n dt_1\cdots dt_p \\ +\sum\limits_{j=1}^{p-1} (-1)^{p-j} \z_n^\star ({\bf m}_j^v)\L_{1,\overleftarrow{({\bf m+1})}_{j+1,p-1}}(1-x)+\z_n^\star(({\bf m}_{p-1},1)^v) \end{array} \right\}.
\end{align*}
Then by the induction hypothesis, we have
\begin{align*}
\z_n^\star(({\bf m}_{p-1},1)^v;x)&=\frac{(-1)^p n}{m_1!\cdots m_{p-1}!} \int\nolimits_{E_p(x)} \Omega_p(1-x)\Omega_{p,p-1}^{m_{p-1}+1}\cdots \Omega_{2,1}^{m_1+1}t_1^n dt_1\cdots dt_p \\ &\quad+\sum\limits_{j=1}^{p-1} (-1)^{p-j} \z_n^\star ({\bf m}_j^v)\L_{1,\overleftarrow{({\bf m+1})}_{j+1,p-1}}(1-x)+\z_n^\star(({\bf m}_{p-1},1)^v).
\end{align*}
Further, if $m_p\geq 1$, we assume the formula (\ref{2.17}) holds for index $(m_1,m_2+1,\ldots,m_{p-1}+1,m_p)^v$. By (\ref{2.16}), it is easy to obtain that
\begin{align*}
\frac{d}{dx} \z_n^\star({\bf m}_{p}^v;x)&=\frac{d}{dx}\left\{\begin{array}{l} \frac{(-1)^p n}{m_1!\cdots m_{p}!} \int\nolimits_{E_p(x)} \Omega_p^{m_p+1}(1-x)\Omega_{p,p-1}^{m_{p-1}+1}\cdots \Omega_{2,1}^{m_1+1}t_1^n dt_1\cdots dt_p \\ \quad+\sum\limits_{j=1}^{p} (-1)^{p-j} \z_n^\star ({\bf m}_j^v)\L_{\overleftarrow{({\bf m+1})}_{j+1,p}}(1-x) \end{array} \right\}\\
&=\frac{\z_n^\star(({\bf m}_{p-1},m_{p})^v)-\z_n^\star(({\bf m}_{p-1},m_{p})^v;x)}{1-x}.
\end{align*}
Hence, the desired evaluation is obtained.\hfill$\square$

\section{Formulas of K-Y MZVs} \label{sec3}
In this section we will prove several explicit formulas of (single-parametric) K-Y MZVs and multiple polylogarithms, and find some explicit relations among MZVs and MZSVs.
\subsection{General Results}

\begin{thm}\label{thm3.1} For ${\bf m}=(m_1,\ldots,m_p)\in\N_0^p$ with $m_1\geq 1$ and index ${\bf k}=(k_1,\ldots,k_r)$,
\begin{align}\label{3.1}
&\sum\limits_{n=1}^\infty \frac{\zeta_{n-1}(k_2,\ldots,k_r)\zeta^\star_n ({\bf m}_p^v;x)}{n^{k_1+1}}\nonumber\\
&=(-1)^{p} \sum\limits_{{\bf k'},j\geq 0}c_{\bf k}({\bf k'};j)\sum\limits_{j_0+\cdots+j_p=j,\atop j_0,\ldots,j_p \geq 0} \frac{(-1)^{j_0}}{j_0!}\left\{\prod\limits_{l=1}^p \binom{m_l+j_l}{j_l}\right\}\nonumber\\&\quad\quad\quad\quad\quad\quad\quad\quad\quad\quad\times\log^{j_0}(1-x){\rm Li}_{\overleftarrow{\left({\bf m+j+1}\right)}_{1,p},{\bf k'}}(1-x)\nonumber\\
&\quad + \sum\limits_{j=1}^{p} (-1)^{p-j} \zeta({\bf k}\circledast (1,{\bf m}_j^v)^\star ) {\rm Li}_{\overleftarrow{\left({\bf m+1}\right)}_{1,j+1}}(1-x),
\end{align}
where $x\in(0,1)$, ${\bf j}:=(j_1,\ldots,j_p)$, ${\bf k'}$ and $c_{\bf k}({\bf k'};j)$ were defined in (\ref{2.11}).
\end{thm}
\pf Multiplying (\ref{2.17}) by $\frac{\zeta_{n-1}(k_2,\ldots,k_r)}{n^{k_1+1}}$ and summing with respect to $n$, we have
\begin{align*}
&\sum\limits_{n=1}^\infty \frac{\zeta_{n-1}(k_2,\ldots,k_r)\zeta^\star_n ({\bf m}_p^v;x)}{n^{k_1+1}}\nonumber\\
&=\frac{(-1)^{p}}{m_1!\cdots m_p!} \int\nolimits_{E_p(x)} \Omega_p^{m_p+1}(1-x)\Omega_{p,p-1}^{m_{p-1}+1}\cdots \Omega_{2,1}^{m_1+1}{\rm Li}_{\bf k}(t_1)  dt_1\cdots dt_p\nonumber\\
&\quad + \sum\limits_{j=1}^{p} (-1)^{p-j} \sum\limits_{n=1}^\infty \frac{\zeta_{n-1}(k_2,\ldots,k_r)\zeta^\star_n ({\bf m}_j^v)}{n^{k_1+1}}{\rm Li}_{m_p+1,m_{p-1}+1,\ldots,m_{j+1}+1}(1-x).
\end{align*}
Then with the help of formula (\ref{2.12}), we may easily deduce the desired evaluation by an elementary calculation. \hfill$\square$

\begin{cor}\label{cor3.2} For positive integers $k,r$ and ${\bf m}=(m_1,\ldots,m_p)\in\N_0^p$ with $m_1\geq 1$,
\begin{align}\label{3.2}
&\sum\limits_{n=1}^\infty \frac{\zeta_{n-1}(\{1\}_{r-1})\zeta^\star_n ({\bf m}_p^v;x)}{n^{k+1}}\nonumber\\
&=(-1)^{p+k-1}\sum\limits_{j_1+\cdots+j_k=r,\atop j_1,\ldots,j_k\geq 0}\sum\limits_{i_0+\cdots+i_p=j_k,\atop i_0,\ldots,i_p \geq 0} \frac{(-1)^{i_0}}{i_0!}\left\{\prod\limits_{l=1}^p \binom{m_l+i_l}{i_l}\right\}\log^{i_0}(1-x)\nonumber\\&\quad\quad\quad\quad\quad\quad\quad\quad\quad\quad\quad\quad\times{\rm Li}_{\overleftarrow{\left({\bf m+i+1}\right)}_{1,p},\overrightarrow{({\bf j+1})}_{1,k-1}}(1-x)\nonumber\\
&\quad+(-1)^p\sum\limits_{j=0}^{k-2} (-1)^j \zeta(k-j,\{1\}_{r-1}){\rm Li}_{\overleftarrow{\left({\bf m+1}\right)}_{1,p},\{1\}_j}(1-x)\nonumber\\
&\quad + \sum\limits_{j=1}^{p} (-1)^{p-j} \zeta((k,\{1\}_{r-1})\circledast(1,{\bf m}_j^v)^\star){\rm Li}_{\overleftarrow{\left({\bf m+1}\right)}_{j+1,p}}(1-x),
\end{align}
where $x\in(0,1)$, ${\bf i}:=(i_1,\ldots,i_p)$ and ${\bf j}:=(j_1,\ldots,j_{k-1})$.
\end{cor}
\pf Corollary \ref{cor3.2} follows immediately from Theorem \ref{thm3.1} by setting $k_1=k,k_2=\cdots=k_r=1$ with the the help of (\ref{2.15}).\hfill$\square$

\begin{cor}\label{cor-c1} For integers $a\geq 1, b\geq 0$ and ${\bf m}=(m_1,\ldots,m_p)\in\N_0^p$ with $m_1\geq 1$,
\begin{align}\label{ec1}
&\sum\limits_{n=1}^\infty \frac{\zeta_{n-1}(\{1\}_{a-1},2,\{1\}_b)\zeta^\star_n ({\bf m}_p^v;x)}{n^{2}}\nonumber\\
&=(-1)^p\sum_{j=0}^a (-1)^j \binom{j+b+1}{j}\z(j+b+2)\sum_{i_0+\cdots+i_p=a-j,\atop i_0,\ldots,i_p\geq 0} \frac{(-1)^{i_0}}{i_0!}\left\{\prod\limits_{l=1}^p \binom{m_l+i_l}{i_l}\right\}\nonumber\\
&\quad\quad\quad\quad\quad\quad\quad\quad\quad\quad\quad\quad\quad\quad\quad\quad\quad\quad\times\log^{i_0}(1-x){\rm Li}_{\overleftarrow{({\bf m+i+1})}_{1,p}}(1-x)\nonumber\\
&-(-1)^{a+p} \sum_{j=0}^{b+1} \binom{j+a}{j}\sum_{i_0+\cdots+i_p=b+1-j,\atop i_0,\ldots,i_p\geq 0} \frac{(-1)^{i_0}}{i_0!}\left\{\prod\limits_{l=1}^p \binom{m_l+i_l}{i_l}\right\}\nonumber\\
&\quad\quad\quad\quad\quad\quad\quad\quad\quad\quad\quad\quad\quad\quad\quad\quad\times\log^{i_0}(1-x){\rm Li}_{\overleftarrow{({\bf m+i+1})}_{1,p},a+1+j}(1-x)\nonumber\\
&\quad+\sum_{j=1}^p (-1)^{p-j}\z((\{1\}_{a},2,\{1\}_b)\circledast(1,{\bf m}_j^v)^\star){\rm Li}_{\overleftarrow{({\bf m+1})}_{j+1,p}}(1-x),
\end{align}
where $x\in(0,1)$, ${\bf i}:=(i_1,\ldots,i_p)$.
\end{cor}
\pf Corollary \ref{cor-c1} follows immediately from Theorem \ref{thm3.1} by setting $k_1=k,k_2=\cdots=k_r=1$ with the the help of (\ref{eb2}).\hfill$\square$

Clearly, from Theorem \ref{thm3.1}, Corollaries \ref{cor3.2} and \ref{cor-c1}, we can find many interesting relations between K-Y MZVs and MZVs. For example, in (\ref{3.1}), (\ref{3.2}) and (\ref{ec1}), setting $x=0$ yield
\begin{align}\label{3.3}
& \sum\limits_{{\bf k'},j\geq 0}c_{\bf k}({\bf k'};j)\sum\limits_{j_1+\cdots+j_p=j,\atop j_0,\ldots,j_p \geq 0} \left\{\prod\limits_{l=1}^p \binom{m_l+j_l}{j_l}\right\}\zeta\left(\overleftarrow{\left({\bf m+j+1}\right)}_{1,p},{\bf k'}\right)\nonumber\\
&=\sum\limits_{j=1}^{p} (-1)^{j+1}\zeta({\bf k}\circledast(1,{\bf m}_j^v)^\star)\zeta\left(\overleftarrow{\left({\bf m+1}\right)}_{j+1,p}\right),
\end{align}
where ${\bf j}=(j_1,\ldots,j_p)$, ${\bf k'}$  and $c_{\bf k}({\bf k'};j)$ were defined in (\ref{2.11}). And
\begin{align}\label{3.4}
&(-1)^{k-1}\sum\limits_{j_1+\cdots+j_k=r,\atop j_1,\ldots,j_k\geq 0}\sum\limits_{i_1+\cdots+i_p=j_k,\atop i_1,\ldots,i_p \geq 0} \left\{\prod\limits_{l=1}^p \binom{m_l+i_l}{i_l}\right\}\zeta\left(\overleftarrow{\left({\bf m+i+1}\right)}_{1,p},\overrightarrow{({\bf j+1})}_{1,k-1}\right)\nonumber\\
&\quad+\sum\limits_{j=0}^{k-2} (-1)^j \zeta(k-j,\{1\}_{r-1})\zeta\left(\overleftarrow{\left({\bf m+1}\right)}_{1,p},\{1\}_j\right)\nonumber\\
&= \sum\limits_{j=1}^{p} (-1)^{j+1} \zeta((k,\{1\}_{r-1})\circledast(1,{\bf m}_j^v)^\star)\zeta\left(\overleftarrow{\left({\bf m+1}\right)}_{j+1,p}\right).
\end{align}
where ${\bf i}:=(i_1,\ldots,i_p)$ and ${\bf j}:=(j_1,\ldots,j_{k-1})$. And
\begin{align}\label{ec2}
&\sum_{j=1}^p (-1)^{j+1}\z((\{1\}_{a},2,\{1\}_b)\circledast(1,{\bf m}_j^v)^\star)\z\left(\overleftarrow{({\bf m+1})}_{j+1,p}\right)\nonumber\\
&=\sum_{j=0}^a (-1)^j \binom{j+b+1}{j}\z(j+b+2)\sum_{i_1+\cdots+i_p=a-j,\atop i_1,\ldots,i_p\geq 0} \left\{\prod\limits_{l=1}^p \binom{m_l+i_l}{i_l}\right\}\z\left(\overleftarrow{({\bf m+i+1})}_{1,p}\right)\nonumber\\
&-(-1)^{a} \sum_{j=0}^{b+1} \binom{j+a}{j}\sum_{i_1+\cdots+i_p=b+1-j,\atop i_1,\ldots,i_p\geq 0} \left\{\prod\limits_{l=1}^p \binom{m_l+i_l}{i_l}\right\}\z\left(\overleftarrow{({\bf m+i+1})}_{1,p},a+1+j\right),
\end{align}
where ${\bf i}:=(i_1,\ldots,i_p)$. In above three formulas (\ref{3.3}), (\ref{3.4}) and (\ref{ec2}), $m_p>1$.

Next, we use the identities (\ref{3.4}) and (\ref{ec2}) to establish some explicit formulas involving MZVs and MZSVs. Setting $r=k=1$ in (\ref{3.4}) gives
\begin{align}\label{3.5}
&\sum\limits_{j=1}^{p} (-1)^{j+1} \zeta(m_p+1,\ldots,m_{j+1}+1) \zeta^\star (2,{\bf m}_j^v)\nonumber\\
&= \sum\limits_{j=1}^p (m_j+1) \zeta(m_p+1,\ldots,m_{j+1}+1,m_j+2,m_{j-1}+1,\ldots,m_1+1).
\end{align}
Setting $a=1$ and $b=0$ in (\ref{ec2}) yields
\begin{align}\label{ec3}
&\sum_{j=1}^p (-1)^{j+1} \zeta(m_p+1,\ldots,m_{j+1}+1)  \left\{\z(2)\z^\star(2,{\bf m}^v_j)-\z^\star(2,2,{\bf m}^v_j)\right\}\nonumber\\
&=\z(2)\sum_{j=1}^p (m_j+1)\z(m_p+1,\ldots,m_{j+1}+1,m_j+2,m_{j-1}+1,\ldots,m_1+1)\nonumber\\
&\quad+\sum_{j=1}^p (m_j+1)\z(m_p+1,\ldots,m_{j+1}+1,m_j+2,m_{j-1}+1,\ldots,m_1+1,2)\nonumber\\
&\quad-2\z(3)\z(m_p+1,\ldots,m_1+1)+2\z(m_p+1,\ldots,m_1+1,3),
\end{align}
where we used the fact that
\[\z((1,2)\circledast(1,{\bf m}^v_j)^\star)=\z(2)\z^\star(2,{\bf m}^v_j)-\z^\star(2,2,{\bf m}^v_j).\]
If letting $p=1$ and $m_1=m\in \N_0$ in (\ref{3.5}) and (\ref{ec3}), by an elementary calculation we obtain
\begin{align*}
\z^\star(2,\{1\}_m)=(m+1)\z(m+2)
\end{align*}
and
\begin{align*}
\z(2)\z^\star(2,\{1\}_m)-\z^\star(2,2,\{1\}_m)&=(m+1)\{\z(2)\z(m+2)+\z(m+2,2)\}\\&\quad\quad-2\z(3,m+1)-2\z(m+4).
\end{align*}
Hence,
\begin{align*}
\z^\star(2,2,\{1\}_m)=2\z(3,m+1)+2\z(m+4)-(m+1)\z(m+2,2).
\end{align*}
Further, from (\ref{3.5}) and (\ref{ec3}) we arrive at
\begin{align}\label{ecd1}
&\sum_{j=1}^p (-1)^{j+1} \zeta(m_p+1,\ldots,m_{j+1}+1)  \z^\star(2,2,{\bf m}^v_j)\nonumber\\
&=2\z(3)\z(m_p+1,\ldots,m_1+1)-2\z(m_p+1,\ldots,m_1+1,3)\nonumber\\
&\quad-\sum_{j=1}^p (m_j+1)\z(m_p+1,\ldots,m_{j+1}+1,m_j+2,m_{j-1}+1,\ldots,m_1+1,2).
\end{align}

Taking $m_1=\cdots=m_p=m\in\N$ in (\ref{3.5}) and noting the fact that $(m,\{m+1\}_{p-1})^v=(\{\{1\}_{m-1},2\}_{p-1},\{1\}_m)$, we obtain
\begin{align*}
&(m+1)\sum\limits_{a+b=p-1,\atop a,b\geq 0} \zeta(\{m+1\}_a,m+2,\{m+1\}_b)\\&=\sum\limits_{j=1}^p (-1)^{j+1}\zeta^\star(\{2,\{1\}_{m-1}\}_j,1)\zeta(\{m+1\}_{p-j}).
\end{align*}
Since the left hand side of above formula can be rewritten as the products of Riemann zeta values, hence, the MZSV $\zeta^\star(\{2,\{1\}_{m-1}\}_j,1)$ on the right hand side can be evaluated by Riemann zeta values. Thus, we have
\begin{align*}
&(m+1)\sum\limits_{a+b=p-1,\atop a,b\geq 0} \zeta(\{m+1\}_a,m+2,\{m+1\}_b)\\
&=\sum\limits_{j=1}^p (-1)^{j+1}\zeta((m+1)j+1)\zeta(\{m+1\}_{p-j}),
\end{align*}
where we used the well-known identity (two different proofs are given in \cite{OW2006,Zl2005}).
\begin{align*}
\zeta^\star(\{2,\{1\}_{m-1}\}_n,1)=(m+1)\zeta((m+1)n+1)\quad (n,m\in \N).
\end{align*}

Letting $m_1=\cdots=m_r=1,m_{r+1}=2,m_{r+2}=\cdots=m_p=1\ (1\leq r \leq p-1)$ in (\ref{3.5}) yields
\begin{align*}
&2\sum\limits_{a+b=r-1,\atop a,b \geq 0} \zeta(\{2\}_{p-r-1},3,\{2\}_a,3,\{2\}_b)\nonumber\\&+2\sum\limits_{a+b=p-r-2,\atop a,b \geq 0} \zeta(\{2\}_{a},3,\{2\}_b,3,\{2\}_r)+3\zeta(\{2\}_{p-r-1},4,\{2\}_r)\nonumber\\
&=\sum\limits_{a+b=r-1,\atop a,b \geq 0}(-1)^a\zeta^\star(\{2\}_{a+1},1)\zeta(\{2\}_{p-r-1},3,\{2\}_b)\nonumber\\
&\quad+(-1)^{r}\sum\limits_{a+b=p-r-1,\atop a,b \geq 0}(-1)^a\zeta^\star(\{2\}_{r+1},1,\{2\}_a,1)\zeta(\{2\}_{b}).
\end{align*}
Note that in \cite{OZ2008}, Ohno and Zudilin gave a brief evaluation ($m_1,m_2\in\N$)
\begin{align*}
\zeta^\star(\{2\}_{m_1},1,\{2\}_{m_2},1)=4\zeta^\star(2m_1+1,2m_2+1)-2\zeta(2m_1+2m_2+2).
\end{align*}
Zagier \cite{DZ2012} found explicit formulas for $\zeta^\star (\{2\}_a,3,\{2\}_b)\ (a,b\in \N_0)$ in terms of rational linear combinations of products $\zeta(m)\pi^{2n}$ with $m+2n=2a+2b+3$. After Zagier's original work, several other proofs have appeared in the literature, see for example \cite{L2012,KP2013}. The formula of $\zeta^\star (\{2\}_a,3,\{2\}_b)$ played an important role in Brown's proof \cite{Bro2012} of the Hoffman conjecture \cite{H1997} that every multiple zeta value is a $\mathbb{Q}$-linear combination of values $\zeta(s_1,\cdots,s_r)$ for which each $s_i$ is either 2 or 3.

Putting $r=1,k=2,m_1=\cdots=m_p=1$ in (\ref{3.4}) gives
\begin{align*}
&\sum\limits_{j=1}^{p} (-1)^{j+1}\zeta^\star(3,\{2\}_{j-1},1)\zeta(\{2\}_{p-j})\\
&=\zeta(2)\zeta(\{2\}_{p})-\zeta(\{2\}_{p+1})-2\sum\limits_{a+b=p-1,\atop a,b\geq 0} \zeta(\{2\}_a,3,\{2\}_b,1).
\end{align*}
Note that the explicit evaluation of MZSVs $\zeta^\star(3,\{2\}_{j-1},1)$ can be found in Zhao \cite{Zp2016}. It should be emphasized that Zhao \cite{Zp2016} also gave more general evaluations for some special MZSVs.

\subsection{Two Special Cases}

In this subsection, we will prove the conclusion that the K-Y MZVs
\[\zeta(\{1\}_r\circledast (1,\{\{1\}_{m-1},2\}_{p-1},\{1\}_m)^\star)\]
and
\[\zeta((\{1\}_m,2,\{1\}_b)\circledast (1,\{\{1\}_{m-1},2\}_{p-1},\{1\}_m)^\star)\]
can be expressed in terms of the products of Riemann zeta values.

 For any index ${\bf m}=(m_1,m_2,\ldots,m_p)$ with $m_1>1$, we define the following \emph{multiple Hurwitz zeta function}
\begin{align*}
\zeta_{HZ} (m_1,m_2,\ldots,m_p;a+1):=\sum\limits_{n_1>n_2>\cdots>n_p>0}\frac{1}{(n_1+a)^{m_1}(n_2+a)^{m_2}\cdots (n_p+a)^{m_p}},
\end{align*}
where $a\neq -1,-2,-3,\ldots.$ For convenience, we set
$\zeta_{HZ} (\emptyset;a+1):=1.$
It is clear that if $a=0$, then
\begin{align*}
\zeta_{HZ} (m_1,m_2,\ldots,m_p;1)=\zeta (m_1,m_2,\ldots,m_p).
\end{align*}
\begin{thm}\label{thm3.3} For $k,m_j\in\N_0\ (j=1,2,\ldots,p-1)$ and $m_p\in\N$,
\begin{align}\label{e3}
&\sum\limits_{i_1+i_2+\cdots+i_p=k,\atop i_1,i_2,\cdots,i_p\geq 0}\left\{\prod\limits_{j=1}^p \binom{m_j+i_j}{i_j}\right\} \zeta\left(\overleftarrow{({\bf m+i+1})}_{1,p}\right)\nonumber\\
&=\frac{(-1)^k}{k!} \lim\limits_{a\rightarrow 0} \frac{\partial^k}{\partial a^k} \left\{\zeta_{HZ}\left(\overleftarrow{({\bf m+1})}_{1,p};a+1\right)\right\},
\end{align}
\end{thm}
\pf We note that the right hand side of formula (\ref{e3}) (not limit) is equal to
\begin{align}\label{e4}
&\frac{(-1)^k}{k!}\frac{\partial^k}{\partial a^k} \left\{\zeta_{HZ}\left(\overleftarrow{({\bf m+1})}_{1,p};a+1\right)\right\}\\
&=\frac{(-1)^k}{k!}\frac{\partial^k}{\partial a^k} \left\{\zeta_{HZ}(m_p+1,\ldots,m_2+1,m_1+1;a+1)\right\}\nonumber\\
&=\frac{(-1)^k}{k!} \sum\limits_{n_1>n_2>\cdots>n_p>0} \frac{\partial^k}{\partial a^k}\left\{ \frac{1}{(n_1+a)^{m_p+1}\cdots(n_{p-1}+a)^{m_2+1}(n_p+a)^{m_1+1}} \right\}\nonumber\\
&=\frac{(-1)^k}{k!} \sum\limits_{n_1>n_2>\cdots>n_p>0} \sum\limits_{i_1+i_2+\cdots+i_p=k\atop i_1,i_2,\ldots,i_p\geq 0} \frac{k!}{i_1!i_2!\cdots i_p!} \prod\limits_{j=1}^p \frac{\partial^{i_j}}{\partial a^{i_j}}\left( \frac{1}{(n_{p+1-j}+a)^{m_j+1}}\right)\nonumber\\
&=\sum\limits_{i_1+i_2+\cdots+i_p=k,\atop i_1,i_2,\cdots,i_p\geq 0}\left\{\prod\limits_{j=1}^p \binom{m_j+i_j}{k_j}\right\} \zeta_{HZ}\left(\overleftarrow{({\bf m+k+1})}_{1,p};a+1\right),
\end{align}
where we used the general Leibniz rule
\begin{align*}
\left(\prod\limits_{j=1}^p f_j\right)^{(k)}=\sum\limits_{k_1+k_2+\cdots+k_p=k,\atop k_1,k_2,\ldots,k_p\in\N_0} \frac{k!}{k_1!k_2!\cdots k_p!} \prod\limits_{j=1}^p (f_j)^{(k_j)}.
\end{align*}
Hence, letting $a\rightarrow 0$ in (\ref{e4}) yields the desired result.\hfill$\square$

\begin{thm}\label{thm3.4} For $k,m,p\in\N$,
\begin{align}\label{e5}
&\frac{(-1)^k}{k!}\frac{\partial^k}{\partial a^k} \left\{\zeta_{HZ}(\{m+1\}_{p};a+1)\right\}\nonumber\\ &=\sum\limits_{c_1+2c_2+\cdots+pc_p=p,\atop c_1,c_2,\ldots,c_p \geq 0} \left\{\prod\limits_{j=1}^p \frac{(-1)^{(j-1)c_j}}{c_j!j^{c_j}} \right\}\nonumber\\ & \quad \quad \quad \quad \quad \times \sum\limits_{k_1+k_2+\cdots+k_{{|\bf c|}_p}=k,\atop k_1,k_2,\ldots,k_{{|\bf c|}_p}\geq 0} \prod_{i=1}^p \prod\limits_{j_i={|\bf c|}_{i-1}+1}^{{|\bf c|}_i} \binom{im+i-1+k_{j_i}}{k_{j_i}} \zeta(im+i+k_{j_i};a+1),
\end{align}
where $\prod\limits_{j=1}^0(\cdot):=1$ and ${|\bf c|}_i:=c_1+c_2+\cdots+c_i,\quad {|\bf c|}_0:=0.$
\end{thm}
\pf In \cite[Eq. (4.7)]{Xu2017}, we gave the recurrence formula
\begin{align}\label{e6}
\zeta_{HZ}(\{m+1\}_p;a+1)=\frac{(-1)^{p-1}}{p}\sum\limits_{i=0}^{p-1} (-1)^i \zeta_{HZ}((p-i)(m+1);a+1)\zeta_{HZ}(\{m+1\}_i;a+1).
\end{align}
According to \cite[Eq. (2.44)]{Riordan58}, the complete Bell polynomials $Y_n(\cdot)$ satisfy the recurrence
\begin{align}\label{e7}
Y_0=1\,,\quad
Y_p(x_1,x_2,\ldots,x_p)=\sum_{j=0}^{p-1}\binom{p-1}{j}x_{p-j}Y_j(x_1,x_2,\ldots,x_j)
    \,,\quad n\geq1\,,
\end{align}
and
\begin{align}\label{e8}
Y_p(x_1,x_2,\ldots,x_p)=\sum\limits_{c_1+2c_2+\cdots+pc_p=p,\atop c_1,c_2,\ldots,c_p \geq 0} \frac{p!}{c_1!c_2!\cdots c_p!} \left(\frac{x_1}{1!}\right)^{c_1}\left(\frac{x_2}{2!}\right)^{c_2}\cdots \left(\frac{x_p}{p!}\right)^{c_p}.
\end{align}
Letting $x_k=(-1)^{k-1}(k-1)!\zeta_{HZ}(k(m+1);a+1)\ (k=1,2,\cdots,p)$ in (\ref{e7}), then comparing it with (\ref{e6}), we obtain
\begin{align}\label{e9}
\zeta_{HZ}(\{m+1\}_p;a+1)=\frac{Y_p\left(0!\zeta_{HZ}(m+1;a+1),\ldots,(-1)^{p-1}(p-1)!\zeta_{HZ}(p(m+1);a+1) \right)}{p!}.
\end{align}
Hence, by (\ref{e8}),
\begin{align}\label{e10}
\zeta_{HZ}(\{m+1\}_p;a+1)=\sum\limits_{c_1+2c_2+\cdots+pc_p=p,\atop c_1,c_2,\ldots,c_p \geq 0} \left\{\prod\limits_{j=1}^p \frac{(-1)^{(j-1)c_j}\zeta_{HZ}^{c_j}(j(m+1);a+1)}{c_j!j^{c_j}} \right\}.
\end{align}
Thus, differentiating (\ref{e10}) $k$ times with respect to $a$ and using general Leibniz rule, we deduce the desired result.\hfill$\square$

In particular, we compute the three cases
\begin{align}
&\frac{(-1)^k}{k!}\frac{\partial^k}{\partial a^k} \left\{\zeta_{HZ}(\{m+1\};a+1)\right\}=\binom{m+k}{k}\zeta_{HZ}(m+k+1;a+1),\label{eq1}\\
&\frac{(-1)^k}{k!}\frac{\partial^k}{\partial a^k} \left\{\zeta_{HZ}(\{m+1\}_2;a+1)\right\}\nonumber\\
&=\frac{1}{2}\sum\limits_{k_1+k_2=k,\atop k_1,k_2\geq 0}\binom{m+k_1}{k_1}\binom{m+k_2}{k_2} \zeta_{HZ}(m+k_1+1;a+1)\zeta_{HZ}(m+k_2+1;a+1)\nonumber\\&\quad-\frac{1}{2}\binom{2m+k+1}{k}\zeta_{HZ}(2m+k+2;a+1),\label{eq2}\\
&\frac{(-1)^k}{k!}\frac{\partial^k}{\partial a^k} \left\{\zeta_{HZ}(\{m+1\}_3;a+1)\right\}\nonumber\\
&=\frac 1{6} \sum\limits_{k_1+k_2+k_3=k,\atop k_1,k_2,k_3\geq 0} \binom{m+k_1}{k_1}\binom{m+k_2}{k_2}\binom{m+k_3}{k_3}\zeta_{HZ}(m+k_1+1;a+1)\nonumber\\
&\quad\quad\quad\quad\quad\quad\quad\quad\quad\quad\times\zeta_{HZ}(m+k_2+1;a+1)\zeta_{HZ}(m+k_3+1;a+1)\nonumber\\
&\quad-\frac{1}{2}\sum\limits_{k_1+k_2=k,\atop k_1,k_2\geq 0} \binom{m+k_1}{k_1}\binom{2m+k_2+1}{k_2} \zeta_{HZ}(m+k_1+1;a+1)\zeta_{HZ}(2m+k_2+2;a+1)\nonumber\\
&\quad+\frac{1}{3}\binom{3m+k+2}{k}\zeta_{HZ}(3m+k+3;a+1).\label{eq3}
\end{align}

Using Theorems \ref{thm3.3} and \ref{thm3.4}, we can prove the following result of a special K-Y MZVs in terms of Riemann zeta values.
\begin{thm} Let $r,m$ and $p$ be positive integers, then
\begin{align}\label{rc1}
\zeta(\{1\}_r\circledast (1,\{\{1\}_{m-1},2\}_{p-1},\{1\}_m)^\star)\in\mathbb{Q}[\text{\rm Riemann Zeta Values}].
\end{align}
\end{thm}
\pf Letting $k=1$ and $m_1=\cdots=m_p=m\ (m\in\N)$ in (\ref{3.4}), and using (\ref{e3}) and (\ref{e5}), we give
\begin{align*}
&\sum\limits_{j=1}^p (-1)^{j+1} \zeta(\{m+1\}_{p-j})\zeta(\{1\}_r\circledast(1,\{\{1\}_{m-1},2\}_{j-1},\{1\}_m)^\star)\\
&=\sum\limits_{i_1+\cdots+i_p=r,\atop i_1,\ldots,i_p \geq 0} \left\{\prod\limits_{l=1}^p \binom{m+i_l}{i_l}\right\}\zeta(m+i_p+1,\ldots,m+i_2+1,m+i_1+1)\\
&\quad\quad\in\mathbb{Q}[\text{\rm Riemann Zeta Values}].
\end{align*}
Then with the help of the fact $\zeta(\{m+1\}_{p})\in\mathbb{Q}[\text{\rm Riemann Zeta Values}]$, the formula (\ref{rc1}) holds. Thus, we complete the proof.\hfill$\square$

\begin{thm} For $m,p\in \N$ and $b\in \N_0$,
\begin{align}\label{ec4}
\zeta((\{1\}_m,2,\{1\}_b)\circledast (1,\{\{1\}_{m-1},2\}_{p-1},\{1\}_m)^\star)\in\mathbb{Q}[\text{\rm Riemann Zeta Values}].
\end{align}
\end{thm}
\pf Setting $a=m_1=m_2=\cdots=m_p=m\in \N$ in (\ref{ec2}), we have
\begin{align*}
&\sum\limits_{j=1}^p (-1)^{j+1} \zeta(\{m+1\}_{p-j})\zeta((\{1\}_m,2,\{1\}_b)\circledast(1,\{\{1\}_{m-1},2\}_{j-1},\{1\}_m)^\star)\\
&=\sum_{j=0}^m (-1)^j\binom{j+b+1}{j}\z(j+b+2)\\ &\quad\quad\quad\quad\quad\quad\quad\times\sum_{j_1+\cdots+j_p=m-j,\atop i_1,\ldots,i_p\geq 0} \left\{\prod\limits_{l=1}^p \binom{m+i_l}{i_l}\right\}\zeta(m+i_p+1,\ldots,m+i_2+1,m+i_1+1)\\
&\quad-(-1)^m \sum_{i_0+i_1+\cdots+i_p=b+1,\atop i_0,i_1,\ldots,i_p\geq 0} \left\{\prod\limits_{l=0}^p \binom{m+i_l}{i_l}\right\}\zeta(m+i_p+1,\ldots,m+i_1+1,m+i_0+1)\\
&\quad\quad\quad\quad\quad\quad\in\mathbb{Q}[\text{\rm Riemann Zeta Values}].
\end{align*}
Then, using (\ref{e3}) and (\ref{e5}), we complete the proof.\hfill$\square$

By elementary calculations, we get three cases:
\begin{align}\label{rc2}
&\zeta(\{1\}_r\circledast(\{1\}_{m},2,\{1\}_m)^\star)\nonumber\\
&=\sum\limits_{n=1}^\infty \frac{\zeta_{n-1}(\{1\}_{r-1})\zeta^\star_n(\{1\}_{m-1},2,\{1\}_m)}{n^2}\nonumber\\
&=\frac{1}{2}\binom{2m+r+1}{r}\zeta(2m+r+2)\nonumber\\&\quad-\frac{1}{2}\sum\limits_{k_1+k_2=r,\atop k_1,k_2\geq 1} \binom{m+k_1}{k_1}\binom{m+k_2}{k_2} \zeta(m+k_1+1)\zeta(m+k_2+1),
\end{align}
\begin{align}\label{rc3}
&\zeta(\{1\}_r\circledast(\{1\}_{m},2,\{1\}_{m-1},2,\{1\}_m)^\star)\nonumber\\
&=\sum\limits_{n=1}^\infty \frac{\zeta_{n-1}(\{1\}_{r-1})\zeta^\star_n(\{1\}_{m-1},2,\{1\}_{m-1},2,\{1\}_m)}{n^2}\nonumber\\
&=\frac{1}{3}\binom{3m+r+2}{r}\zeta(3m+r+3)\nonumber\\
&\quad-\frac{1}{2}\sum\limits_{k_1+k_2=r,\atop k_1,k_2\geq 1} \binom{m+k_1}{k_1}\binom{2m+k_2+1}{k_2} \zeta(m+k_1+1)\zeta(2m+k_2+2)\nonumber\\
&\quad+\frac{1}{6}\sum\limits_{k_1+k_2+k_3=r,\atop k_1,k_2,k_3\geq 1} \binom{m+k_1}{k_1}\binom{m+k_2}{k_2}\binom{m+k_3}{k_3}\nonumber \\&\quad\quad\quad\quad\times\zeta(m+k_1+1)\zeta(m+k_2+1) \zeta(m+k_3+1).
\end{align}
and
\begin{align}\label{ec5}
&\zeta((\{1\}_{m},2,\{1\}_b)\circledast(\{1\}_{m+1})^\star)\nonumber\\
&=\sum\limits_{n=1}^\infty \frac{\zeta_{n-1}((\{1\}_{m-1},2,\{1\}_b)\zeta^\star_n(\{1\}_{m})}{n^2}\nonumber\\
&=\sum_{j=0}^{m-1} (-1)^j \binom{j+b+1}{j}\binom{2m-j}{m}\z(j+b+2)\z(2m+1-j)\nonumber\\
&\quad-\frac{(-1)^m}{2} \sum_{k_1+k_2=b+1,\atop k_1,k_2\geq 1} \binom{m+k_1}{m}\binom{m+k_2}{m} \z(m+k_1+1)\z(m+k_2+1)\nonumber\\
&\quad+\frac{(-1)^m}{2}\binom{2m+b+2}{b+1} \z(2m+b+3).
\end{align}

The following conjecture seems to follow from formula (\ref{e5}), but we have not yet worked it out in detail.
\begin{con}\label{con1} For any $k,m,p\in\N$,
\begin{align}\label{rc4}
&\zeta(\{1\}_r\circledast (1,\{\{1\}_{m-1},2\}_{p-1},\{1\}_m)^\star)\nonumber\\
&=\sum\limits_{n=1}^\infty \frac{\zeta_{n-1}(\{1\}_{r-1})\zeta^\star_n\left(\{\{1\}_{m-1},2\}_{p-1},\{1\}_m\right)}{n^2}\nonumber\\
&=(-1)^{p+1}\sum\limits_{c_1+2c_2+\cdots+pc_p=p,\atop c_1,c_2,\ldots,c_p \geq 0} \left\{\prod\limits_{j=1}^p \frac{(-1)^{(j-1)c_j}}{c_j!j^{c_j}} \right\}\nonumber\\ &\times \sum\limits_{k_1+k_2+\cdots+k_{{|\bf c|}_p}=r,\atop k_1,k_2,\ldots,k_{{|\bf c|}_p}\geq 1} \prod_{i=1}^p \prod\limits_{j_i={|\bf c|}_{i-1}+1}^{{|\bf c|}_i} \binom{im+i-1+k_{j_i}}{k_{j_i}} \zeta(im+i+k_{j_i}).
\end{align}
\end{con}

If setting $r=1$ in Conjecture \ref{con1}, we can get the well-known identity
\begin{align*}
\zeta^\star( \{2,\{1\}_{m-1}\}_{p},1)=(m+1)\zeta(p(m+1)+1 ).
\end{align*}
Moreover, the formulas (\ref{rc2}) and (\ref{rc3}) can also be obtained by the Conjecture \ref{con1} with $p=2$ and $3$, respectively.

\section{Ap\'{e}ry-Type Multiple Zeta Star Values}\label{sec4}

In this section, we use a similar method in the above section to evaluate the following Ap\'{e}ry-type multiple zeta star values $\z_B^\star({\bf k})$ which involve center binomial coefficient,
\begin{align}
\zeta^\star_B({\bf k})\equiv\zeta^\star_B(k_1,k_2,\ldots,k_r):=\sum\limits_{n=1}^\infty \frac{\zeta_n^\star(k_2,\ldots,k_r)}{n^{k_1}4^n}\binom{2n}{n},
\end{align}
where $k_1,k_2,\ldots,k_r$ are positive integers. In particular, if $r=1$ and $k_1=k$ then
\begin{align*}
\zeta^\star_B({k}):=\sum\limits_{n=1}^\infty \frac{\binom{2n}{n}}{n^{k}4^n}.
\end{align*}

For any index ${\bf m}=(m_1,\ldots,m_p)$, we define the classical multiple polylogarithm function with $r$-complex variables by
\begin{align*}
{\rm Li}_{\bf m}({\bf x})\equiv{\L}_{{{m_1},{m_2}, \cdots ,{m_p}}}( x_1,x_2,\ldots,x_p): = \sum\limits_{n_1>n_2>\cdots>n_p\geq 1} {\frac{{{x_1^{{n_1}}}}x_2^{n_2}\cdots x_p^{n_p}}{{n_1^{{m_1}}n_2^{{m_2}} \cdots n_p^{{m_p}}}}}.
\end{align*}
where ${\bf x}=(x_1,\ldots,x_p)$ with$|x_1\cdots x_j|\leq 1\ (1\leq j\leq p)$ and $(x_1,m_1)\neq (1,1)$.

It is clear that if $x_j=\pm 1$ then ${\rm Li}_{\bf m}({\bf x})$ becomes to (alternating) MZVs. According to definition, we deduce
\begin{align}\label{4.2}
{\rm Li}_{\bf m}({\bf x})=\frac{1}{m_1!\cdots m_p!} \int\nolimits_{E_p(0)} \prod\limits_{j=1}^p \frac{\log^{m_j-1}\left(\frac{t_{j-1}}{t_j}\right)}{(x_1\cdots x_j)^{-1}-t_j}dt_1\cdots dt_p,
\end{align}
where $t_0:=1$.

Now, we prove two formulas of iterated integrals.
\begin{thm} For nonnegative integers $m_1,\ldots,m_{p-1}$ and $m_p>0$,
\begin{align}\label{a6}
&\int\nolimits_{E_p(0)} \Omega_p^{m_p+1}(1)\Omega_{p,p-1}^{m_{p-1}+1}\cdots \Omega_{2,1}^{m_1+1}\log\left(\frac{2}{1+\sqrt{1-t_1}}\right) dt_1\cdots dt_p\nonumber\\
&= m_1!\cdots m_p!2^{|{\bf m}|_p}\log(2) \sum\limits_{\sigma_j \in \{\pm 1\} \atop j=1,2,\ldots,p } \L_{\overleftarrow{({\bf m+1})}_{1,p}}\left(\sigma_p, \cat_{i=1}^{p-1}\{\sigma_{p+1-i} \sigma_{p-i}\} \right)\nonumber\\
&\quad+ m_1!\cdots m_p!2^{|{\bf m}|_p} \sum\limits_{\sigma_j \in \{\pm 1\} \atop j=1,2,\ldots,p } \L_{\overleftarrow{({\bf m+1})}_{1,p},1}\left(\sigma_p, \cat_{i=1}^{p-1}\{\sigma_{p+1-i} \sigma_{p-i}\},-\sigma_1 \right),
\end{align}
where $\cat_{i=j}^k\{a_i\}$ abbreviates the concatenated argument sequence $a_{j},a_{j+1},\ldots,a_{k}$, with the convention $\cat_{i=j}^k\{a_i\}:=\emptyset$ for $k<j$.
\end{thm}
\pf Let $A$ denote the left hand side of (\ref{a6}). Applying the change of variables $t_j=1-t_{p+1-j}\ (j=1,2,\ldots,p)$, we have
\begin{align*}
A= \int\nolimits_{E_p(0)} \frac{\log\left(\frac{1}{t_1}\right)^{m_p}\log\left(\frac{t_1}{t_2}\right)^{m_{p-1}}\cdots \log\left(\frac{t_{p-1}}{t_p}\right)^{m_1}\log\left(\frac{2}{1+\sqrt{t_p}}\right)}{(1-t_1)(1-t_2)\cdots (1-t_p)}dt_1dt_2\cdots dt_p.
\end{align*}
Then, letting $t_j\rightarrow t_j^2\ (j=1,2,\ldots,p)$ yields
\begin{align*}
A&=2^{p+|{\bf m}|_p} \int\nolimits_{E_p(0)} t_1\cdots t_p\frac{\log\left(\frac{1}{t_1}\right)^{m_p}\log\left(\frac{t_1}{t_2}\right)^{m_{p-1}}\cdots \log\left(\frac{t_{p-1}}{t_p}\right)^{m_1}\log\left(\frac{2}{1+t_p}\right)}{(1-t_1^2)(1-t_2^2)\cdots (1-t_p^2)}dt_1\cdots dt_p\\
&=2^{p+|{\bf m}|_p}\log(2)\int\nolimits_{E_p(0)} t_1\cdots t_p\frac{\log\left(\frac{1}{t_1}\right)^{m_p}\log\left(\frac{t_1}{t_2}\right)^{m_{p-1}}\cdots \log\left(\frac{t_{p-1}}{t_p}\right)^{m_1}}{(1-t_1^2)(1-t_2^2)\cdots (1-t_p^2)}dt_1\cdots dt_p\\
&\quad-2^{p+|{\bf m}|_p}\int\nolimits_{E_p(0)} t_1\cdots t_p\frac{\log\left(\frac{1}{t_1}\right)^{m_p}\log\left(\frac{t_1}{t_2}\right)^{m_{p-1}}\cdots \log\left(\frac{t_{p-1}}{t_p}\right)^{m_1}\log\left({1+t_p}\right)}{(1-t_1^2)(1-t_2^2)\cdots (1-t_p^2)}dt_1\cdots dt_p\\
&=2^{|{\bf m}|_p} \log(2) \sum\limits_{\sigma_j\in\{\pm 1\}\atop j=1,2,\ldots,p}\int\nolimits_{E_p(0)}\frac{\log\left(\frac{1}{t_1}\right)^{m_p}\log\left(\frac{t_1}{t_2}\right)^{m_{p-1}}\cdots \log\left(\frac{t_{p-1}}{t_p}\right)^{m_1}}{(\sigma_p^{-1}-t_1)(\sigma_{p-1}^{-1}-t_2)\cdots (\sigma_1^{-1}-t_p)}dt_1\cdots dt_p\\
&\quad-2^{|{\bf m}|_p} \sum\limits_{\sigma_j\in\{\pm 1\}\atop j=1,2,\ldots,p}\int\nolimits_{E_p(0)}\frac{\log\left(\frac{1}{t_1}\right)^{m_p}\log\left(\frac{t_1}{t_2}\right)^{m_{p-1}}\cdots \log\left(\frac{t_{p-1}}{t_p}\right)^{m_1}\log(1+t_p)}{(\sigma_p^{-1}-t_1)(\sigma_{p-1}^{-1}-t_2)\cdots (\sigma_1^{-1}-t_p)}dt_1\cdots dt_p.
\end{align*}
Then, using (\ref{4.2}) and noting that
\[\log(1+t_p)=\int\limits_{0}^{t_{p}} \frac{1}{1+t_{p+1}}dt_{p+1},\]
we obtain the desired evaluation by a direct calculation.\hfill$\square$

\begin{thm} For nonnegative integers $m_1,\ldots,m_{p-1}$ and $m_p>0$,
\begin{align}\label{a7}
&\int\nolimits_{E_p(0)} \Omega_p^{m_p+1}(1)\Omega_{p,p-1}^{m_{p-1}+1}\cdots \Omega_{2,1}^{m_1+1}\frac{\log\left(\frac{1+\sqrt{1-t_1}}{2\sqrt{1-t_1}}\right)}{\sqrt{1-t_1}} dt_1\cdots dt_p\nonumber\\
&=m_1!\cdots m_p!2^{|{\bf m}|_p} \sum\limits_{|{\bf i}|_p=1\atop i_1,i_2,\ldots,i_p\geq 0} \prod\limits_{l=1}^{p} \binom{m_l+i_l}{i_l}\sum\limits_{\sigma_j \in \{\pm 1\} \atop j=1,2,\ldots,p } \sigma_1\L_{\overleftarrow{({\bf m+i+1})}_{1,p}}\left(\sigma_p, \cat_{i=1}^{p-1}\{\sigma_{p+1-i} \sigma_{p-i}\} \right) \nonumber\\
&\quad-m_1!\cdots m_p!2^{|{\bf m}|_p}\log(2) \sum\limits_{\sigma_j \in \{\pm 1\} \atop j=1,2,\ldots,p }\sigma_1 \L_{\overleftarrow{({\bf m+1})}_{1,p}}\left(\sigma_p, \cat_{i=1}^{p-1}\{\sigma_{p+1-i} \sigma_{p-i}\} \right)\nonumber\\
&\quad- m_1!\cdots m_p!2^{|{\bf m}|_p} \sum\limits_{\sigma_j \in \{\pm 1\} \atop j=1,2,\ldots,p }\sigma_1 \L_{\overleftarrow{({\bf m+1})}_{1,p},1}\left(\sigma_p, \cat_{i=1}^{p-1}\{\sigma_{p+1-i} \sigma_{p-i}\},-\sigma_1 \right).
\end{align}
\end{thm}
\pf The proof of (\ref{a7}) is similar as the proof of (\ref{a6}). Let $B$ denote the left hand side of (\ref{a7}). By an elementary calculation, we can find that
\begin{align*}
B&=2^{p+|{\bf m}|_p} \int\nolimits_{E_p(0)} t_1\cdots t_{p-1}\frac{\log\left(\frac{1}{t_1}\right)^{m_p}\log\left(\frac{t_1}{t_2}\right)^{m_{p-1}}\cdots \log\left(\frac{t_{p-1}}{t_p}\right)^{m_1}\log\left(\frac{1+t_p}{2t_p}\right)}{(1-t_1^2)(1-t_2^2)\cdots (1-t_p^2)}dt_1\cdots dt_p\\
&=2^{p+|{\bf m}|_p} \int\nolimits_{E_p(0)} t_1\cdots t_{p-1}\frac{\log\left(\frac{1}{t_1}\right)^{m_p}\log\left(\frac{t_1}{t_2}\right)^{m_{p-1}}\cdots \log\left(\frac{t_{p-1}}{t_p}\right)^{m_1}\log\left(1+t_p\right)}{(1-t_1^2)(1-t_2^2)\cdots (1-t_p^2)}dt_1\cdots dt_p\\
&\quad-2^{p+|{\bf m}|_p}\log(2) \int\nolimits_{E_p(0)} t_1\cdots t_{p-1}\frac{\log\left(\frac{1}{t_1}\right)^{m_p}\log\left(\frac{t_1}{t_2}\right)^{m_{p-1}}\cdots \log\left(\frac{t_{p-1}}{t_p}\right)^{m_1}}{(1-t_1^2)(1-t_2^2)\cdots (1-t_p^2)}dt_1\cdots dt_p\\
&\quad+2^{p+|{\bf m}|_p} \int\nolimits_{E_p(0)} t_1\cdots t_{p-1}\frac{\log\left(\frac{1}{t_1}\right)^{m_p}\log\left(\frac{t_1}{t_2}\right)^{m_{p-1}}\cdots \log\left(\frac{t_{p-1}}{t_p}\right)^{m_1}\log\left(\frac{1}{t_p}\right)}{(1-t_1^2)(1-t_2^2)\cdots (1-t_p^2)}dt_1\cdots dt_p.
\end{align*}
We note the fact that
\begin{align*}
\log\left(\frac{1}{t_p}\right)&=\log\left(\frac{t_{p-1}}{t_p}\cdots \frac{t_1}{t_2}\frac{1}{t_1}\right)=\sum_{j=1}^p \log\left(\frac{t_{j-1}}{t_j} \right)\quad\quad (t_0:=1)\\
&=\sum\limits_{i_1+i_2+\cdots+i_p=1\atop i_1,i_2,\ldots,i_p\geq 0}\log\left(\frac{1}{t_1}\right)^{i_p}\log\left(\frac{t_1}{t_2}\right)^{i_{p-1}}\cdots \log\left(\frac{t_{p-1}}{t_p}\right)^{i_1}.
\end{align*}
Hence, with the help of formula (\ref{4.2}) we may easily deduce the desired result.\hfill$\square$

Next, we prove two recurrence relations for Ap\'{e}ry-type multiple zeta values $\zeta^\star_B(\cdots)$. First, by integration by parts, we get
\begin{align*}
&\frac{(-1)^{p+r} n^{r+1}}{m_1!\cdots m_p!} \int\nolimits_{E_{p+r}(x)} \Omega_{p+r}^{m_p+1}(1-x)\Omega_{p+r,p+r-1}^{m_{p-1}+1}\cdots \Omega_{r+2,r+1}^{m_1+1}\Omega_{r+1,r}\cdots \Omega_{2,1}t_1^n \prod\limits_{j=1}^{p+r}dt_j\\
&=\frac{(-1)^p n}{m_1!\cdots m_p!} \int\nolimits_{E_p(x)} \Omega_p^{m_p+1}(1-x)\Omega_{p,p-1}^{m_{p-1}+1}\cdots \Omega_{2,1}^{m_1+1}t_1^n dt_1\cdots dt_p\\&\quad-\sum\limits_{j=1}^r n^{r+1-j} (-1)^{p+r-j} \L_{\overleftarrow{({\bf m+1})}_{1,p},\{1\}_{r-j}}(1-x).
\end{align*}
Hence, the (\ref{2.17}) can be rewritten in the form
\begin{align}\label{4.5}
&\z^\star_n({\bf m}_p^v;x )\nonumber\\
&=\frac{(-1)^{p+r} n^{r+1}}{m_1!\cdots m_p!} \int\nolimits_{E_{p+r}(x)} \Omega_{p+r}^{m_p+1}(1-x)\Omega_{p+r,p+r-1}^{m_{p-1}+1}\cdots \Omega_{r+2,r+1}^{m_1+1}\Omega_{r+1,r}\cdots \Omega_{2,1}t_1^n \prod\limits_{j=1}^{p+r}dt_j\nonumber\\
&\quad+\sum\limits_{j=1}^r n^{r+1-j} (-1)^{p+r-j} \L_{\overleftarrow{({\bf m+1})}_{1,p},\{1\}_{r-j}}(1-x)\nonumber\\
&\quad+ \sum\limits_{j=1}^{p} (-1)^{p-j}\z^\star_n({\bf m}_j^v)\L_{\overleftarrow{({\bf m+1})}_{j+1,p}}(1-x).
\end{align}
\begin{thm}\label{thm4.3} For integers $m_2,\ldots,m_{p-1},r\geq 0$ and $m_1,m_p\geq 1$,
\begin{align}\label{a8}
&\sum\limits_{j=1}^p (-1)^{j-1} \zeta\left(\overleftarrow{({\bf m+1})}_{j+1,p}\right) \z^\star_B(r+2,{\bf m}_j^v)\nonumber\\
&=\sum\limits_{j=1}^r (-1)^{r-j}\zeta\left(\overleftarrow{(\bf m+1)}_{1,p},\{1\}_{r-j}\right)\z^\star_B(j+1)\nonumber\\
&\quad+(-1)^r2^{|{\bf m}|_p+1}\log(2) \sum\limits_{\sigma_j \in \{\pm 1\} \atop j=1,2,\ldots,p+r } \L_{\overleftarrow{({\bf m+1})}_{1,p},\{1\}_{r}}\left(\sigma_{p+r}, \cat_{i=1}^{p+r-1}\{\sigma_{p+r+1-i} \sigma_{p+r-i}\} \right)\nonumber\\
&\quad+2^{|{\bf m}|_p+1} \sum\limits_{\sigma_j \in \{\pm 1\} \atop j=1,2,\ldots,p+r } \L_{\overleftarrow{({\bf m+1})}_{1,p},\{1\}_{r+1}}\left(\sigma_{p+r}, \cat_{i=1}^{p+r-1}\{\sigma_{p+r+1-i} \sigma_{p+r-i}\},-\sigma_1 \right).
\end{align}
\end{thm}
\pf Setting $x=0$ in (\ref{4.5}) gives
\begin{align}\label{a9}
&\sum\limits_{j=1}^{p} (-1)^{j-1}\z^\star_n({\bf m}_j^v)\zeta\left(\overleftarrow{({\bf m+1})}_{j+1,p}\right)\nonumber\\
&=\frac{(-1)^{r} n^{r+1}}{m_1!\cdots m_p!} \int\nolimits_{E_{p+r}(0)} \Omega_{p+r}^{m_p+1}(1)\Omega_{p+r,p+r-1}^{m_{p-1}+1}\cdots \Omega_{r+2,r+1}^{m_1+1}\Omega_{r+1,r}\cdots \Omega_{2,1}t_1^n dt_1\cdots dt_{p+r}\nonumber\\
&\quad+\sum\limits_{j=1}^r n^{r+1-j} (-1)^{r-j} \zeta\left(\overleftarrow{(\bf m+1)}_{1,p},\{1\}_{r-j}\right).
\end{align}
Multiplying it by $\frac{\binom{2n}{n}}{n^{r+2}4^n}$ and summing with respect to $n$, and noting that the well-known formula (see \cite{ChenH16,Leh1985})
\begin{equation*}
\sum_{n=1}^\infty\frac{\binom{2n}{n}}{4^n n}t^n
    =2\log\left(\frac{2}{1+\sqrt{1-t}}\right),
\end{equation*}
we obtain
\begin{align*}
&\sum\limits_{j=1}^p (-1)^{j-1} \zeta\left(\overleftarrow{({\bf m+1})}_{j+1,p}\right) \sum\limits_{n=1}^\infty \frac{\z^\star_n ({\bf m}_j^v)}{n^{r+2}4^n} \binom{2n}{n}\\
&=\frac{2(-1)^{r} }{\prod\limits_{j=1}^p m_j!} \int\nolimits_{E_{p+r}(0)} \Omega_{p+r}^{m_p+1}(1)\Omega_{p+r,p+r-1}^{m_{p-1}+1}\cdots \Omega_{r+2,r+1}^{m_1+1}\Omega_{r+1,r}\cdots \Omega_{2,1}\log\left(\frac{2}{1+\sqrt{1-t_1}}\right)dt_1\cdots dt_{p+r}\nonumber\\
&\quad+\sum\limits_{j=1}^r (-1)^{r-j}\zeta\left(\overleftarrow{(\bf m+1)}_{1,p},\{1\}_{r-j}\right)\sum\limits_{n=1}^\infty \frac{\binom{2n}{n}}{n^{j+1}4^n}.
\end{align*}
Then by (\ref{a6}), the desired formula can be established directly.\hfill$\square$

\begin{thm}\label{thm4.4}  For integers $m_2,\ldots,m_{p-1},r\geq 0$ and $m_1,m_p\geq 1$,
\begin{align}\label{a10}
&\sum\limits_{j=1}^p (-1)^{j-1} \zeta\left(\overleftarrow{({\bf m+1})}_{j+1,p}\right) \sum\limits_{n=1}^\infty \frac{\z^\star_n ({\bf m}_j^v)H_n}{n^{r+1}4^n} \binom{2n}{n}\nonumber\\
&=\sum\limits_{j=1}^r (-1)^{r-j}\zeta\left(\overleftarrow{(\bf m+1)}_{1,p},\{1\}_{r-j}\right)\z^\star_B (j,1)\nonumber\\
&\quad+(-1)^r 2^{|{\bf m}|_p+1} \sum\limits_{|{\bf i}|_r+|{\bf j}|_p=1\atop i_1,\ldots,i_r,j_1,\ldots,j_p\geq 0} \prod\limits_{l=1}^{p} \binom{m_l+j_l}{j_l}\nonumber\\&\quad\quad\quad\quad\times\sum\limits_{\sigma_j \in \{\pm 1\} \atop j=1,2,\ldots,p+r} \sigma_1\L_{\overleftarrow{({\bf m+j+1})}_{1,p},\overleftarrow{({\bf i+1})}_{1,r}}\left(\sigma_{p+r}, \cat_{i=1}^{p+r-1}\{\sigma_{p+r+1-i} \sigma_{p+r-i}\} \right)\nonumber\\
&\quad-(-1)^r2^{|{\bf m}|_p+1}\log(2) \sum\limits_{\sigma_j \in \{\pm 1\} \atop j=1,2,\ldots,p+r }\sigma_1 \L_{\overleftarrow{({\bf m+1})}_{1,p},\{1\}_{r}}\left(\sigma_{p+r}, \cat_{i=1}^{p+r-1}\{\sigma_{p+r+1-i} \sigma_{p+r-i}\} \right)\nonumber\\
&\quad-(-1)^r2^{|{\bf m}|_p+1} \sum\limits_{\sigma_j \in \{\pm 1\} \atop j=1,2,\ldots,p+r }\sigma_1 \L_{\overleftarrow{({\bf m+1})}_{1,p},\{1\}_{r+1}}\left(\sigma_{p+r}, \cat_{i=1}^{p+r-1}\{\sigma_{p+r+1-i} \sigma_{p+r-i}\},-\sigma_1 \right),
\end{align}
where $H_n=\zeta_n(1)$ stands for the classical harmonic number.
\end{thm}
\pf The proof of (\ref{a10}) is similar as the proof of (\ref{a8}). Multiplying (\ref{a9}) by $\frac{H_n\binom{2n}{n}}{n^{r+1}4^n}$ and summing with respect to $n$, and noting that the formula (see Chen \cite{ChenH16})
\begin{equation*}
\sum_{n=1}^\infty\frac{H_n\binom{2n}{n}}{4^n}t^n
    =\frac{2}{\sqrt{1-t}}\log\left(\frac{1+\sqrt{1-t}}{2\sqrt{1-t}}\right),
\end{equation*}
with the help of (\ref{a7}), we prove the desired result.\hfill$\square$

From (\ref{a8}) and (\ref{a10}), we can get the following corollary.
\begin{cor} For an index ${\bf m}=(m_1,\ldots,m_p)$ and positive integer $r$,
\begin{align*}
\zeta_B^\star(r,{\bf m}_p^v)\in\mathbb{Q}[\text{\rm Alternating MZVs with depth}\leq p+r].
\end{align*}
\end{cor}

From Theorems \ref{thm4.3} and \ref{thm4.4}, by direct calculations, we deduce two specific examples:
\begin{align*}
&\z^\star_B(2,2,1)=\sum_{n=1}^\infty \frac{\zeta^\star_n(2,1)}{n^24^n}\binom{2n}{n}\\
&=2\zeta^2(2)\log(2)+4\zeta(2)(\zeta(\bar 2,1)+\zeta(2,\bar 1))\\
&\quad-8\log(2)\left( \zeta(2,2)+\zeta(2,\bar 2)+\zeta(\bar 2,\bar 2)+\zeta(\bar 2,2)\right)\\
&\quad-8\left( \zeta(2,2,\bar 1)+\zeta(2,\bar 2,1)+\zeta(\bar 2,\bar 2,\bar 1)+\zeta(\bar 2,2,1)\right)\\
&=\frac{75}{8}\zeta(5)-4\zeta(4)\log(2)-3\zeta(2)\zeta(3)
\end{align*}
and
\begin{align*}
&\sum_{n=1}^\infty \frac{H^2_n}{n^24^n}\binom{2n}{n}=2\z^\star_B(2,1,1)-\z^\star_B(2,2)\\
&=\zeta(2)\z^\star_B(1,1)-8(\zeta(3,1)-\zeta(3,\bar 1)+\zeta(\bar 3,\bar 1)-\zeta(\bar 3,1))\\
&\quad-4(\zeta(2,2)-\zeta(2,\bar 2)+\zeta(\bar 2, \bar 2)-\zeta(\bar 2,2))\\
&\quad+4\log(2)(\zeta(2,1)-\zeta(2,\bar 1)+\zeta(\bar 2,\bar 1)-\zeta(\bar 2,1))\\
&\quad+4(\zeta(2,1,\bar 1)-\zeta(2,\bar 1,1)+\zeta(\bar 2,\bar 1,\bar 1)-\zeta(\bar 2,1,1))\\
&=32{\rm Li}_4(1/2)-14\zeta(4)+7\zeta(3)\log(2)-8\zeta(2)\log^2(2)+\frac{4}{3}\log^4(2),
\end{align*}
where we used Au's Mathematica package [1,
Appedix A] containing the explicit expressions of all values of alternating MZVs with weight $\leq 8$. The last series was also found in \cite[Exa. 2.4]{WX2019}.

We end the paper by the following conjecture based on our computations.
\begin{con} For ${\bf k}=(k_1,\ldots,k_r)\in \N^r$, the {\rm MZBSVs} $\zeta^\star_B({\bf k})$ can be expressed in terms of alternating {\rm MZVs}.

\end{con}

{\bf Acknowledgments.} The author expresses his deep gratitude to Professors Masanobu Kaneko, Jianqiang Zhao, Weiping Wang and Nobuo Sato for
their valuable comments and encouragement. This work was supported by the China Scholarship Council (No. 201806310063) and the Scientific Research Foundation for Scholars of Anhui Normal University.

 \end{document}